\documentclass{article}

\usepackage{amsfonts}
\usepackage{amssymb}
\usepackage{epsfig}
\usepackage{amsthm}

\newcommand{\N}{{\mathbb N}}
\newcommand{\Z}{{\mathbb Z}}
\newcommand{\X}{{\cal X}}
\newcommand{\G}{{\hat G}}
\renewcommand{\H}{{\hat H}}
\newcommand{\ga}{$\Gamma(G,{\cal A})$ }
\renewcommand{\O}{{\cal O}}
\newcommand{\R}{{\cal R}}

\newcommand{\A}{{\hat A}}

\newtheorem{lemma}{Lemma}[section]
\newtheorem{cor}{Corollary}
\newtheorem{thm}{Theorem}

\newcounter{remarks}
\newcommand{\remark}{\par \refstepcounter{remarks}{\noindent \bf Remark \arabic{remarks}. }}

\begin{document}
\begin{center}
{\large ON RESIDUALIZING HOMOMORPHISMS PRESERVING QUASICONVEXITY

\vspace{.5cm}
 Ashot Minasyan}

\vspace{.3cm}
Department of Mathematics

Vanderbilt University

Nashville, TN 37240, USA

aminasyan@gmail.com
\end{center}

\begin{abstract}
$H$ is called a $G$-subgroup of a hyperbolic group $G$ if for any finite subset $M\subset G$ there exists
 a homomorphism
from $G$ onto a non-elementary hyperbolic group $G_1$ that is surjective on $H$ and injective on $M$.
In his paper in $1993$ A. Ol'shanskii
gave a description of all $G$-subgroups in any given non-elementary hyperbolic group $G$.
Here we show that for the same class of $G$-subgroups the finiteness assumption on $M$
(under certain natural conditions) can be replaced by an assumption of quasiconvexity.
\end{abstract}

\section{Introduction}
Suppose $G$ is a group with some finite generating set $\cal A$ (assume that $\cal A$ is {\it symmetrized}, i.e.
${\cal A}^{-1}={\cal A}$). Then one can construct the Cayley graph
\ga which is a geodesic metric space. Assume $\eta \ge 0$ is some number.
A subset $Q$ of $G$ (or of \ga) is said to be $\eta$-quasiconvex (or just quasiconvex) if for any pair of elements
$u,v \in Q$ and any geodesic segment $p$ connecting $u$ and $v$, $p$ belongs to a closed $\eta$-neighborhood
of the subset $Q$ in \ga.

If, in addition, the group $G$ is (word) hyperbolic, the fact that a subset of this group is quasiconvex
does not depend on the choice of a generating set $\cal A$ (see \cite{Gromov}).

If $S \subseteq G$,
$C_H(S)$ will denote the centralizer subgroup of $S$ in $H$, i.e. $$C_H(S)=\{h\in H~|~hg=gh~~\forall~ g \in S\}~.$$
If $A,B$ are subsets of the group $G$, their product is the subset of $G$ defined by
$$A \cdot B=AB \stackrel{def}{=} \{ab~|~a\in A,b\in B\}~.$$ $A^{-1} \subset G$ will denote the subset
$\{a^{-1}~|~a \in A\}$.

Now let $G$ be a $\delta$-hyperbolic group, $\delta \ge 0$. Every element $g \in G$ of infinite order
belongs to a unique {\it maximal elementary subgroup} $E(g)$ (a group is called {\it elementary} if it has a cyclic
subgroup of finite index). Then by \cite[Lemmas 1.16 and 1.17]{Olsh2}
$$E(g)=\bigl\{~x \in G~|~xg^nx^{-1}=g^{\pm n}~\mbox{for some $n\in \N$}~\bigr\} ~\mbox{  and  }$$
\begin{equation} E(g)=\bigl\{~x \in G~|~xg^kx^{-1}=g^{l}~\mbox{for some $k,l \in \Z\backslash \{0\}$}~\bigr\}
\label{elemdef}
\end{equation}
It is easy to see that the subgroup $$E^+(g)=\bigl\{~x \in G~|~xg^nx^{-1}=g^{ n}~\mbox{for some $n\in \N$}~\bigr\}$$
is of index at most $2$ in $E(g)$.

For any subgroup $H$ of $G$ denote by $H^0$ the set of elements of infinite order in $H$.
Using the terminology from \cite{Olsh2}, a non-elementary subgroup $H$ of $G$ will be called
a {\it $G$-subgroup} if for any finite subset $M\subset G$ there is a non-elementary hyperbolic quotient $G_1$
of the group $G$, such that the natural homomorphism $G\to G_1$ is surjective on $H$ and injective on $M$.

Denote $E(H)=\bigcap_{x \in H^0} E(x)$. If $H$ is a non-elementary subgroup of $G$, then $E(H)$ is the unique
maximal finite subgroup of $G$ normalized by $H$ (\cite[Prop. 1]{Olsh2}). Hence $H$ acts on $E(H)$ by conjugation
and we have a  homomorphism of $H$ into the permutation group on the set of elements of $E(H)$.
The kernel of that homomorphism is $C_H\bigl(E(H)\bigr)$ which sometimes will be  denoted by $K(H)$.
The index $|H:K(H)|$ is finite because of the finiteness of $E(H)$.

The following characterization of all $G$-subgroups was given in \cite[Thm. 1]{Olsh2}:
\begin{lemma} A non-elementary subgroup $H$ of a hyperbolic group $G$ is a $G$-subgroup if and only if
$E(H)=E(G)$ and $|H:K(H)|=|G:K(G)|$ (i.e. the actions by conjugation of $H$ and $G$ on $E(H)=E(G)$ are similar:
for every $g \in G$ there exists an element $h \in H$ with $gag^{-1}=hah^{-1}$ for all $a \in E(G)$).
\end{lemma}

\vspace{.15cm}
\underline{\bf Definition.} Let $H$ be a subgroup of the group $G$ and $Q\subseteq G$ be a quasiconvex subset.
The subset $Q$ will be called {\it small relatively to} $H$ if for any two finite subsets
$P_1,P_2$ of the group $G$ one has
$$H \nsubseteq \left( P_1\cdot Q^{-1} \cdot Q \cdot P_2\right)~.  \eqno (*)$$

\vspace{.15cm}
As we know, any generating set induces a left-invariant metric on the set of elements of a group.
So, if $G_1$ is a quotient of $G$, the group $G_1$ will be generated by the image of $\cal A$ under the natural
homomorphism $\phi: G \to G_1$. Therefore, later $G_1$ will be assigned the metric corresponding to the generating
set $\phi(\cal A)$.

The main result of this paper is
\begin{thm} \label{mainthm} Let $H_1$, $H_2$, \dots, $H_k$ be $G$-subgroups of a non-elementary hyperbolic group
$G$ and $H_1',\dots,H'_{k'}$ be some non-elementary subgroups of $G$. Assume $Q \subseteq G$ is
an $\eta$-quasiconvex subset (for some $\eta \ge 0$) that is small relatively to $H_i$ for every $i=1,2\dots,k$.
%
Then there exist a group $G_1$ and an epimorphism $\phi: G \to G_1$ such that

$\bf  1)$ $G_1$ is a non-elementary hyperbolic group;

$\bf 2)$ The homomorphism $\phi$ is an isometry between $Q$ and $\phi(Q)$ (if the metrics on $G$ and $G_1$ are
chosen as explained above) and for any quasiconvex subset $S \subseteq Q$, its image $\phi(S)$ is
quasiconvex in $G_1$. In particular, $\phi$ is injective on $Q$;

$\bf 3)$ $\phi$ is surjective on each of the subgroups $H_1,\dots,H_k$, i.e. $\phi(H_i)=G_1$ for each
$i=1,2,\dots,k$;

$\bf 4)$ $\phi$-images of two elements from $Q$ are conjugate in $G_1$ if and only if these elements
are conjugate in $G$;

$\bf 5)$ The centralizer $C_{G_1}\bigl(\phi(a)\bigr)$ for every $a \in Q$ is the $\phi$-image of the
centralizer $C_G(a)$;

$\bf 6)$ $ker \phi$ is a torsion-free subgroup;

$\bf 7)$ $\phi$ induces a bijective map on sets of conjugacy classes of elements having finite orders
in $G$ and $G_1$ respectively;

$\bf 8)$ $\phi(H_1'),\dots,\phi(H'_{k'})$ are non-elementary subgroups of $G_1$;

$\bf 9)$ $E(G_1)=\phi\bigl(E(G)\bigr)$.
\end{thm}

Parts 1)-8) of Theorem \ref{mainthm} were proved by A. Ol'shanskii in \cite[Thm. 2]{Olsh2} for the case when the
subset $Q$ is finite; part 9) was proved in \cite[Thm. 4]{Olsh2} with additional conditions imposed on $G$.
If $card(Q)<\infty$, since each $H_i$ is infinite, the condition $(*)$ becomes trivial and can be omitted
($Q$ will be small relatively to any infinite subgroup).
In our proof  the tools and techniques developed in the paper \cite{Olsh2} will be crucial, so
an interested reader is referred to that article in advance.

\vspace{.15cm}
\underline{\bf Definition.}  A subgroup $H$ of a group $G$ will be called a {\it quasiretract} of $G$ if there
exists a normal subgroup $N \lhd G$ such that $|G:HN|<\infty$ and the intersection $H \cap N$ is finite.

\vspace{.15cm}
In particular, any retract is a quasiretract. The proof of the following lemma is not difficult and will be given
in the beginning of section \ref{proofofthm2}.

\begin{lemma} \label{quasi-quasi} Assume that $G$ is a hyperbolic group and $H$ is a quasiretract of $G$.
Then the subgroup $H$ is quasiconvex in $G$.
\end{lemma}

One can observe that if the group $G$ is torsion-free then any its non-elemen\-tary subgroup will be a $G$-subgroup.
However, by far, not every subgroup in $G$ will be quasiconvex (or a quasiretract).
In the next theorem we show that demanding $Q$ to be small relatively to $H_i$
is necessary if one doesn't impose additional limitations on the subgroups $H_i$, $i=1,2,\dots,k$.

\begin{thm} \label{ness} Let $H$ be an infinite subgroup of a hyperbolic group $G$ and $Q \subset G$ be a
subset (not necessarily quasiconvex). Suppose that $$H \subseteq P_1Q^{-1}QP_2$$ for some finite subsets
$P_1$, $P_2$ of $G$
and there is a group $G_1$ and an epimorphism $\phi: G \to G_1$ such that $\phi$ is surjective on
$H$ (i.e. $\phi(H)=G_1$) and $\phi|_Q$ is a quasiisometry  between $Q$ and $\phi(Q)$.
Then the subgroup $H$ is a quasiretract of $G$.
\end{thm}

In section \ref{*-cond} we investigate some properties of condition $(*)$; we show that if two quasiconvex subsets
satisfy this condition then so do their union and product (lemma \ref{smallunion}). If, in addition, one demands that
$Q^{-1}$ is quasiconvex then we are able to simplify $(*)$: the set $Q^{-1}Q$ in it can be substituted by just $Q$
(corollary \ref{q-1q=q}).

As examples of quasiconvex sets $Q$ satisfying $(*)$ we can consider special
subsets that are finite unions of products of quasiconvex subgroups. More precisely,
suppose $F_1,\dots,F_n$ are quasiconvex subgroups of a hyperbolic group $G$ and
$g_0,g_1,\dots,g_n \in G$. Following \cite{hyp}, the subset
$$P=g_0F_1\cdot \dots \cdot g_{n-1}F_ng_n=\{g_0f_1\cdot \dots \cdot g_{n-1}f_ng_n~|~f_i \in F_i,i=1,\dots,n\}
\subseteq G$$
will be called a {\it quasiconvex product}. The quasiconvex subgroups $F_i$, $i=1,\dots,n$, are said to be
{\it members} of the product $P$. A finite union of quasiconvex products is always a quasiconvex subset
in a hyperbolic group (\cite[Prop. 1,Lemma 2.1]{hyp}).

Let $U=\bigcup_{k=1}^N P_k$ be a finite union of quasiconvex products $P_k$, $k=1,\dots,N$.
A subgroup $F \le G$ will be called a { \it member} of $U$, by definition, if $F$ is a member of $P_k$
for some $1 \le k \le N$. For any such set $U$ we fix its representation as a finite union of quasiconvex
products and fix its members.

\begin{thm} \label{qcex} Suppose $G$ is a hyperbolic group, $H$ is its subgroup and $Q$ is a finite
union of quasiconvex products in $G$ with members $F_1,\dots,F_l$. Assume also that
$$|H:(H\cap gF_jg^{-1})|=\infty~,~\mbox{ for every } g \in G \mbox{ and } j=1,2,\dots,l~.$$
Then for arbitrary two finite subsets
$P_1,P_2 \subset G$, the subgroup $H$ is not contained inside of $P_1\cdot Q^{-1} \cdot Q \cdot P_2~.$
\end{thm}

Let $\cal A$ be an alphabet and $\R_i$ -- subsets of words in ${\cal A}^{\pm 1}$, $i \in \N$, satisfying
$\R_i \subset \R_{i+1}$ for all $i$. Let the groups $G_i$ have presentations
$$G_i=\langle {\cal A}~\|~\R_i\rangle,~ i \in \N.$$
Then $G_{i+1} \cong G_i/N_i$ where $N_i \lhd G_i$ is the normal closure of $\R_{i+1}\backslash\R_i$ in $G_i$, i.e.
there is an epimorphism $\phi_i: G_i \to G_{i+1}$ with $ker(\phi_i)=N_i$.

Set $\R = \bigcup_{i=1}^{\infty} \R_i$.
The group $M$ defined by the presentation
$$M=\langle {\cal A}~\|~\R\rangle$$ is said to be an {\it inductive limit} of the groups $G_i$, $i \in \N$.

Thus we obtain an infinite sequence of epimorphisms
$$G_1 \stackrel{\phi_1}{\to} G_2 \stackrel{\phi_2}{\to} G_3 \stackrel{\phi_3}{\to} \dots$$
and $\displaystyle M=\lim_{\to} (G_i,\phi_i)$.

The statement below follows from theorem \ref{mainthm} and can be compared with \cite[Theorems 1,3]{Olsh-SQ}.

\begin{cor} \label{simpleinductive} Suppose $G$, $H$ are hyperbolic groups and $G$ is non-elementary. Then $H$ can be isomorphically
embedded into some simple quotient $M$ of the group $G$. Moreover, the group $M$ is an inductive limit of
hyperbolic groups.
\end{cor}

D. Osin noticed that using a similar argument one can obtain

\begin{cor} \label{monster} There exists a simple group $M$ that is
a quotient of every non-elementary hyperbolic group and contains every hyperbolic group (isomorphically embedded).
\end{cor}

Recall that a non-trivial proper subgroup $H$ of a group $G$ is called {\it malnormal} if for any
$g \in G \backslash H$ the intersection $H \cap gHg^{-1}$ is trivial.

In the torsion-free case we can strengthen the claim of corollary \ref{simpleinductive} to obtain a
"thrifty" embedding (cf. \cite{Olsh-thrifty}):

\begin{cor} \label{thrifty} Suppose $G$, $H$ are torsion-free hyperbolic groups, $G$ is non-elementary
and $H$ is non-trivial. Then there exists a simple torsion-free quotient $M$ of $G$ and an injective homomorphism
$\pi: H \to M$ such that $\pi(H)$ is malnormal in $M$ and any proper subgroup of $M$ is
conjugate (in $M$) to a subgroup of $\pi(H)$.
\end{cor}

\vspace{.3cm}
{\noindent \large \bf Acknowledgements}

The author is grateful to Professor A. Yu. Ol'shanskii
for suggesting the problem and helpful remarks.

\section{Preliminary Information}
Assume $\bigl(\X,d(\cdot,\cdot)\bigr)$ is a proper geodesic metric space.
If $Q \subset \X$, $N \ge 0$, the closed $N$-neighborhood of $Q$ will be denoted by
$${\cal O}_N (Q) \stackrel{def}{=} \{x\in \X~|~d(x,Q) \le N \}~.$$

If $x,y,w \in \X$, then the number
$$(x|y)_w \stackrel{def}{=} \frac12 \Bigl(d(x,w)+d(y,w)-d(x,y) \Bigr)$$
is called the {\it Gromov product} of $x$ and $y$ with respect to $w$.

Let $abc$ be a geodesic triangle in the space $\X$ and $[a,b]$, $[b,c]$, $[a,c]$ be its sides between the corresponding
vertices. There exist "special" points $O_a \in [b,c]$,\\ $O_b \in [a,c]$, $O_c \in [a,b]$
with the properties:
$d(a,O_b) = d(a,O_c) = \alpha$, $d(b,O_a) = $ $=d(b,O_c) = \beta$, $d(c,O_a) = d(c,O_b) = \gamma$. From a corresponding
system of linear equations one can find that $\alpha = (b|c)_a$, $\beta = (a|c)_b$, $\gamma = (a|b)_c$. Two points
$O \in [a,b]$ and $O' \in [a,c]$ are called $a$-{\it equidistant} if $d(a,O) = d(a,O') \le \alpha$.
The triangle $abc$ is said to be $\delta$-{\it thin} if for any two points $O,O'$ lying on its sides and
equidistant from one of its vertices, $d(O,O') \le \delta$ holds.

A geodesic $n$-gon in the space $\X$ is said to be $\delta$-{\it slim} if each of its sides belongs to a
closed $\delta$-neighborhood of the union of the others.

Later we will use three equivalent definitions of  {\it hyperbolicity} of the space $\X$ (see \cite{Ghys},\cite{Mihalik}): \\
$1^{\circ}.$ (M. Gromov) There exists $\delta \ge 0$ such that for any four points $x,y,z,w \in \X$ their Gromov products satisfy
$$(x|y)_w \ge min\{(x|z)_w,(y|z)_w\} - \delta~;$$
$2^{\circ}.$ All triangles in $\X$ are {\cal $\delta$-thin} for some $\delta \ge 0$;\\
$3^{\circ}.$ (E. Rips) All triangles in $\X$ are $\delta$-slim for some $\delta \ge 0$.

\vspace{.15cm}
\remark It is easy to see that the definition $3^\circ$ implies that any geodesic $n$-gon in the space $\X$
is $(n-2)\delta$-slim if $n \ge 3$.
\vspace{.15cm}

Now, suppose $G$ is a finitely generated group with a fixed finite generating set $\cal A$. For every element
$g\in G$ its length $|g|_G$ is the length of a shortest word in the alphabet $\cal A$ representing $g$ in $G$.
Then, if $g,h \in G$ we define the distance $d(g,h)=|g^{-1}h|_G$. This distance function can be extended
to a metric on the Cayley graph invariant under the action of $G$ by left multiplication:

$$~\forall~x,y \in \mbox{\ga and } g \in G~~d(g\circ x,g\circ y)=d(x,y)~.$$

This implies that for any elements $x,y,z,w \in G$
$$(x|y)_w=(zx|zy)_{zw}~.$$
Thus the Cayley graph \ga becomes a proper geodesic metric space. There is a natural (metric space) embedding of the
group $G$ into its Cayley graph. Later we will identify $G$ with its image under it.

The group $G$ is called {\it hyperbolic} if \ga is a hyperbolic metric space. This definition does not depend on
the choice of the finite generating set $\cal A$ in $G$ (\cite{Gromov},\cite[5.2.14]{Ghys}),
thus hyperbolicity is a group-theoretical property.
Among well-known examples of hyperbolic groups are free groups of finite rank, finite groups, fundamental groups
of negatively curved compact manifolds, etc.

Further on we will assume that \ga meets $1^{\circ},2^{\circ}$ and $3^{\circ}$ for a fixed (sufficiently large)
$\delta \ge 0$.

For any two points $x,y \in \Gamma(G,{\cal A})$ we fix a geodesic path between them and denote it by $[x,y]$.
Let $p$ be a path in the Cayley graph of $G$. Then $p_-$, $p_{+}$  will denote the startpoint and
the endpoint of $p$, $||p||$ -- its length; $lab(p)$, as usual, will mean the word in the alphabet $\cal A$
written on $p$. $elem(p) \in G$ will denote the element of the group $G$ represented by the word $lab(p)$.
$p^{-1}$ will be the inverse path to $p$, i.e. the path with the same set of points but traced in the opposite direction.

A path $q$ is called $(\lambda,c)$-{ \it quasigeodesic} if there exist $0<\lambda \le 1$, $c \ge 0$, such that
for any subpath $p$ of $q$ the inequality $\lambda ||p|| - c \le d(p_-,p_+)$ holds.\\
In a hyperbolic space quasigeodesics and geodesics with same ends are mutually close~:

\begin{lemma} \label{close} {\normalfont (\cite[5.6,5.11]{Ghys},\cite[3.3]{Mihalik})}
There is a constant $\nu=\nu(\delta,\lambda,c)$ such that for any
$(\lambda,c)$-quasigeodesic path $p$ in $\Gamma(G,{\cal A})$ and a geodesic
$q$ with $p_- = q_-$, $p_+ = q_+$, one has $p \subset {\cal O}_\nu(q)$ and $q \subset
{\cal O}_\nu(p)$. \end{lemma}

An important property of cyclic subgroups in a hyperbolic group states

\begin{lemma} \label{power} {\normalfont (\cite[8.21]{Ghys},\cite[3.2]{Mihalik})}
For any word $w$ representing an element $g \in G$ of infinite order  there
exist constants $\lambda >0$, $c \ge 0$, such that any path with a label $w^m$ in the
Cayley graph of $G$ is $(\lambda,c)$-quasigeodesic for arbitrary integer $m$.
\end{lemma}

It is easy to see that any finite subset of a group is $D/2$-quasiconvex where $D$ is the diameter of that subset.

It follows from lemmas \ref{close} and \ref{power} that any cyclic subgroup of a
hyperbolic group is quasiconvex.

\begin{lemma} \label{q-cunion,product} {\normalfont (\cite[Prop. 3.14]{Zeph},\cite[Lemma 2.1,Prop. 1]{hyp})}
Let $G$ be a hyperbolic group and let $A,B$ be its quasiconvex subsets. Then the subsets $A\cup B$
and $A \cdot B$ are also quasiconvex.
\end{lemma}

If $X_1,X_2,\dots, X_n$ are points in \ga, the notation $X_1X_2 \dots X_n$ will be used for the geodesic $n$-gon
with vertices $X_i$, $i=1,\dots,n$, and sides $[X_i,X_{i+1}]$, $i=1,2,\dots ,n-1$, $[X_n,X_0]$. $[X_1,X_2,\dots,X_n]$
will denote the broken line with these vertices in the corresponding order (i.e. the path $[X_1,X_2,\dots,X_n]$
will consist of consecutively concatenated geodesic segments $[X_i,X_{i+1}]$, $i=1,2,\dots,n$).

\begin{lemma} \label{brokenlines1} {\normalfont (\cite[Lemma 21]{Olsh1})} Let $p = [X_0,X_1,\dots,X_n]$ be a broken
line in $\Gamma(G,{\cal A})$ such that
$||[X_{i-1},X_i]|| \ge C_1$ $\forall~i=1,\dots,n$, and $(X_{i-1}|X_{i+1})_{X_i} \le C_0$ $\forall~i=1,\dots,n-1$,
where $C_0 \ge 14 \delta$, $C_1 > 12(C_0 + \delta)$. Then $p$ is contained in the closed $2C_0$-neighborhood
${\cal O}_{2C_0}([X_0,X_n])$ of the geodesic segment $[X_0,X_n]$.
\end{lemma}

\begin{lemma} \label{brokenlines2} {\normalfont (\cite[Lemma 3.5]{paper2})} In the conditions of lemma
\ref{brokenlines1} the inequality $\|[X_0,X_n]\| \ge \|p\|/2$ holds.
\end{lemma}

\begin{lemma} \label{quadrangle}  {\normalfont (\cite[Lemma 4.1]{paper2})} Consider  a geodesic
quadrangle $X_1X_2X_3X_4$  in \ga with
$d(X_2,X_3)>d(X_1,X_2)+d(X_3,X_4)$. Then there are points \\$U,V \in [X_2,X_3]$ such that
$d(X_2,U) \le d(X_1,X_2)$, $d(V,X_3) \le d(X_3,X_4)$ and  the geodesic subsegment
$[U,V]$ of $[X_2,X_3]$ lies $2\delta$-close to the side $[X_1,X_4]$.
\end{lemma}

\begin{lemma} \label{quasintersection} {\normalfont (\cite[Lemma 4.5]{paper2})} Let $G$ be a $\delta$-hyperbolic
and let $H_i$ be $\varepsilon_i$-quasiconvex subgroups of $G$, $i=1,2$.
If one has $$\sup \{(h_1|h_2)_{1_G}~:~h_1 \in H_1,h_2 \in H_2 \} = \infty$$ then
$card(H_1\cap H_2)=\infty$.
\end{lemma}

\begin{lemma} \label{conj} {\normalfont (\cite[Lemma 4.2]{paper2})} Let $A$ be an infinite $\varepsilon$-quasiconvex
set in $G$  and $g \in G$.
Then if the intersection $A \cap gAg^{-1}$ is infinite, there exists an element $r\in G$ with
$|r|_G \le 4\delta+ 2\varepsilon + 2\varkappa$ such that $g \in ArA^{-1}$, where $\varkappa$ is the length
of a shortest element from $A$.
\end{lemma}

\begin{lemma} \label{newdouble} {\normalfont (\cite[Lemma 4.3]{paper2})}
Let $G$ be a $\delta$-hyperbolic group, $H$ and $K$ -- its subgroups
with $K$ quasiconvex. If ~~$\displaystyle H \subset \bigcup_{j=1}^N Ks_jK$~~ for some $s_1,\dots,s_N \in G$ then
$|H:(K \cap H)| < \infty$ .
\end{lemma}

\begin{lemma} \label{non-periodic} {\normalfont (\cite[8.3.36]{Ghys})}
Any infinite subgroup of a hyperbolic group contains an element of infinite order.
\end{lemma}

If $g \in G^0$, $T(g)$ will be used to denote the set of elements of finite order in the subgroup $E(g)$.

\vspace{.15cm}
{\bf \underline{Definition.}} Let $G$ be a hyperbolic group and $H$ be its non-elementary subgroup.
An element $g \in H^0$ will be called $H$-{\it suitable} if $E(H)=T(g)$ and
$$E(g)=E^+(g)=C_G(g)=T(g)\times\langle g \rangle~.$$

\vspace{.15cm}
In  particular, if the element $g$ is $H$-suitable then $g \in C_H\bigl(E(H)\bigr)$.

Two elements $g,h \in G$ of infinite order are called {\it commensurable} if \\ $g^k=ah^la^{-1}$ for some non-zero
integers $k,l$ and some $a \in G$.

Now we recall the statement of \cite[Lemma 3.8]{Olsh2}:
\begin{lemma} \label{lemma3.8} Every non-elementary subgroup $H$ of a hyperbolic group $G$ contains an infinite set
of pairwise non-commensurable $H$-{\it suitable} elements. 
\end{lemma}

We will need the following modification of \cite[Lemma 3.7]{Olsh2}:
\begin{lemma} \label{lemma3.7} Let $g$ be an $H$-suitable element in a non-elementary subgroup $H$ of
a hyperbolic group $G$. Suppose $l \in \N$ and $K$ is a non-elementary subgroup of $H$.
Then for any number $C_1 \ge 0$ there exist elements  $x_i \in K$, $i=1,\dots,l$, satisfying the following properties:

{\bf 0)} $|x_i|_G>C_1$ for every $i=1,2\dots,l$;

{\bf 1)} $x_i \notin E(g)$ for every $i=1,2\dots,l$;

{\bf 2)} $x_i \in C_G(E(H))$ for every $i=1,2\dots,l$;

{\bf 3)} $ax_i=x_ib$ for $a,b \in E(g)$ implies that $a=b \in E(H)$, $i=1,2\dots,l$;

{\bf 4)} if $a,b \in E(g)$ and $ax_i=x_jb$ for some $i,j \in \{1,\dots,l\}$ then $i=j$.
\end{lemma}

\underline{Proof.}
Indeed, it is shown in the proof of \cite[Lemma 3.7]{Olsh2} that if the elements $g,h_1,\dots,h_l \in H$ are pairwise
non-commensurable in $G$ then for any sufficiently large $t \in \N$, the elements $x_i=h_i^t$ satisfy the conditions
$1)-4)$. By lemma \ref{lemma3.8} we can choose such $h_1,\dots,h_l$ inside of $K$, thus $x_i=h_i^t \in K$.
Obviously, if $t\in \N$ is sufficiently large the property $0)$ will also be satisfied. $\square$

\vspace{.15cm}
In this paper we will also use the concept of {\it Gromov boundary} of a hyperbolic group $G$ (for a detailed theory
the reader is referred to the corresponding chapters in \cite{Ghys},\cite{Bridson} or \cite{Mihalik}).)
In order to define it, call a sequence $(x_i)_{i\in \N} \subset \Gamma(G,{\cal A})$ {\it converging to infinity} if
$$\lim_{i,j \to \infty} (x_i|x_j)_{1_G} = \infty~\mbox{ ($1_G$ is the identity element of the group $G$}).$$

Two sequences $(x_i)_{i\in \N},(y_j)_{i\in \N}$ converging to infinity are said to be equivalent
if $$\lim_{i \to \infty} (x_i|y_i)_{1_G} = \infty~.$$
The points of the boundary $\partial G$ are identified with the equivalence classes of sequences
converging to infinity. (It is easy to see that this definition does not depend on the choice
of a basepoint: instead of $1_G$ one could use any fixed point $p$ of \ga \cite{Mihalik}.)
If $\alpha$ is the equivalence class of $(x_i)_{i\in \N}$ we will write
$\displaystyle \lim_{i \to \infty} x_i = \alpha$.

The space $\partial G$ can be topologized so that it becomes compact, Hausdorff
and metrizable (see \cite{Mihalik},\cite{Ghys}).

Every isometry $\psi$ of the space \ga induces a homeomorphism of $\partial G$ in a natural way:
for every equivalence class of sequences converging to infinity $[(x_i)_{i\in \N}] \in \partial G$
choose a representative $(x_i)_{i\in \N} \subset \Gamma(G,{\cal A})$
and set $$\psi \left(\left[(x_i)_{i\in \N}\right]\right)\stackrel{def}{=}
\left[\left(\psi(x_i)\right)_{i\in \N}\right] \in \partial G~.$$

Left multiplication by elements of the group induces an isometric action of $G$ on \ga. Hence,
$G$ acts homeomorphically on the boundary $\partial G$ as described above.

If $g \in G$ is an element of infinite order then the sequences $(g^i)_{i\in \N}$ and
$(g^{-i})_{i\in \N}$ converge to infinity and we will use the notation
$$\lim_{i \to \infty} g^i = g^\infty,~\lim_{i \to \infty} g^{-i} = g^{-\infty}~.$$

For a subset $A \subseteq \Gamma(G,{\cal A})$ the {\it limit set} $\Lambda (A)$ of $A$ is the collection of the points
$\alpha \in \partial G$ that are limits of sequences (converging to infinity) from $A$.

The following properties of limit sets are well-known:

\begin{lemma} \label{limitsets} {\normalfont (\cite{K-S},\cite{Swenson})}
Suppose $A,B$ are arbitrary subsets of \ga, $g \in G$. Then

(a) $\Lambda (A) = \emptyset$ if and only if $A$ is finite;

(b) $\Lambda (A)$ is a closed subset of the boundary $\partial G$;

(c) $\Lambda(A \cup B) = \Lambda (A) \cup \Lambda (B)$;

(d) $\Lambda(Ag) = \Lambda(A)$, $g \circ \Lambda(A)=\Lambda(gA)$;

\end{lemma}

\begin{lemma} \label{boundedproduct} Suppose $A$ and $B$  are subsets of the hyperbolic group $G$ and
$\Lambda(A)\cap \Lambda(B)=\emptyset$. Then $\displaystyle \sup_{a \in A,~b\in B}\{(a|b)_{1_G}\}< \infty~.$
\end{lemma}

\underline{Proof.} This statement is an easy consequence of the definition of a limit set. Indeed,
assume, by the contrary, that there are sequences of elements $(a_i)_{i\in \N} \subset A$ and $(b_i)_{i\in \N} \subset B$
such that ~$\lim_{i\to \infty} (a_i|b_i)_{1_G} = \infty$. Then the subsets $\{a_i~|~i\in \N\} \subset G$ and
$\{b_i~|~i\in \N\} \subset G$ are infinite, hence each of them has at least one limit point (by lemma
\ref{limitsets}.(a)). Thus, there are subsequences $(a_{i_j})_{j\in \N}$ of $(a_i)$ and $(b_{i_j})_{j\in \N}$ of $(b_i)$
satisfying
$$\lim_{j \to \infty} a_{i_j} = \alpha \in \Lambda(A),~\lim_{j \to \infty} b_{i_j} = \beta \in \Lambda(B)~.$$
But $\lim_{j\to \infty} (a_{i_j}|b_{i_j})_{1_G} = \infty$ by our assumption, hence $\alpha=\beta$. A contradiction
with the condition $\Lambda(A)\cap \Lambda(B)=\emptyset$. $\square$


If $H$ is a subgroup of $G$, it is known that $\Lambda H$ is either empty (if $H$ is finite) or
consists of two distinct points (if $H$ is infinite elementary), or is uncountable (if $H$ is
non-elementary) (\cite{K-S},\cite{Ghys}). In the second case, when there exists $g \in H^0$ such that
$|H:\langle g \rangle| < \infty$, $\Lambda H = \{g^{\infty},g^{-\infty} \}$.

As the hyperbolic group $G$ acts on its boundary, for every subset $\Omega \subset \partial G$ one can
define the stabilizer subgroup by $St_G(\Omega) = \{g \in G~|~g \circ \Omega = \Omega\}$.
For our convenience, we set $St_G(\emptyset) = G$.

It is proved in \cite[thm. 8.3.30]{Ghys} that for any point $\alpha \in \partial G$~~ $St_G(\{\alpha\})$ is
an elementary subgroup of the group $G$ (in fact, if $\alpha=g^\infty$ for some element
of infinite order $g \in G$ then $$St_G(\{g^\infty\})=E^+(g)~,$$
otherwise the subgroup $St_G(\{\alpha\})$ is finite).
$St_G(\{g^\infty,g^{-\infty}\}) = E(g)$.

If $A \subseteq G$ and $\Omega \subseteq \partial G$ then the orbit of $\Omega$ under the action of $A$ will be
denoted $$A \circ \Omega \stackrel{def}{=} \{g \circ \alpha~|~g \in A,\alpha \in \Omega\} \subseteq \partial G~.$$

\begin{lemma} \label{invariant} {\normalfont (\cite[Lemma 3.3]{K-S})} If $H$ is an infinite subgroup of $G$ then
$\Lambda(H)$ contains at least
two distinct points and the set $\Lambda(H)$ is $H$-invariant, i.e.
for every $h \in H$ ~~$h\circ \Lambda(H)=\Lambda(H)$.
\end{lemma}

Thus every infinite  subgroup $H$ of $G$ acts on its limit set $\Lambda(H)$ and this action is
induced by the action of $G$ on the boundary $\partial G$ described above.

\begin{lemma} \label{limitofnormal} {\normalfont (\cite[Lemma 3.8]{K-S})} Let $A$ be an infinite normal subgroup of a
subgroup $H$ in $G$. Then $\Lambda(A)=\Lambda(H)$.
\end{lemma}

For an arbitrary subset $\Omega \subseteq \partial G$ denote by $cl(\Omega)$ its closure in the topology of
the boundary $\partial G$.

\begin{lemma} \label{stabilizer} {\normalfont (\cite[Lemma 6.4]{paper2})} Suppose $\Omega \subset \partial G$
is a  subset having at least two distinct points. Then $\Lambda\bigl(St_G(\Omega)\bigr) \subseteq cl(\Omega)$.
\end{lemma}

Combining lemmas $8.1.(5)$ and $8.2$ from \cite{paper2} together one obtains
\begin{lemma} \label{tame-qc} Assume $H$ is a subgroup of $G$ and $A \subseteq G$ is a non-empty quasiconvex subset. If
$\Lambda(H) \subseteq \Lambda(A)$ then there is a finite subset $P \subset G$ such that
$H \subseteq A \cdot P=\bigcup_{x \in P}Ax$.
\end{lemma}

\section{Proofs of Theorems \ref{ness},\ref{qcex}} \label{proofofthm2}
Suppose $H= \langle {\cal B} \rangle$ is a subgroup of $G$ with a finite generating set $\cal B$.
If $h \in H$, then by $|h|_G$ and $|h|_H$ we will denote the lengths of the element $h$ in the alphabets $\cal A$ and
$\cal B$ respectively. Denote $\N_0=\N \cup \{0\}=\{0,1,2,\dots\}$.
The {\it distortion function} (see \cite{Gromov2}) $D_H: ~\N_0 \to \N_0$ of $H$ in $G$ is defined by
$$D_H(n) = max\{|h|_H~|~h \in H, ~|h|_G \le n \}~.$$
If ${\alpha},{\beta}~: \N_0 \to \N_0$ are two functions then we write ${\alpha} \preceq {\beta}$ if there are constants
$K_1,K_2>0$~: ${\alpha}(n) \le K_1{\beta}(K_2n)$ for every $n \in \N$. $\alpha$ and $\beta$ are said to be
equivalent if ${\alpha} \preceq {\beta}$ and ${\beta} \preceq {\alpha}$.

Evidently, the function $D_H$ does not depend (up to this equivalence) on the choice of finite generating sets
$\cal A$ of $G$ and $\cal B$ of $H$. One can also notice that $D_H(n)$ is always at least linear
(provided that $H$ is infinite).
Fix a linear function $L:\N_0\to\N_0$ (for example, $L(n)=n$). If $D_H \preceq L$,
 $H$ is called {\it undistorted}.

\begin{lemma} \label{quasi-undist} {\normalfont \cite[Lemma 1.6]{hyp}} Let $H$ be a finitely generated subgroup of a
hyperbolic group $G$. Then $H$ is quasiconvex if and only if $H$ is undistorted in $G$.
\end{lemma}

\underline{Proof of lemma \ref{quasi-quasi}.} Let $H$ be a quasiretract of the hyperbolic group $G$
and let $N \lhd G$ satisfy $|G:HN|<\infty$, $card(H\cap N)<\infty$. Denote $\G=HN \le G$.
Then for the quotient group $K=\G/N$ there is a natural epimorphism $\varphi: \G \to K$ such that $K=\varphi(H)$ and
$\displaystyle M \stackrel{def}{=}ker(\varphi) \cap H$ is finite.

The group $\G$ is generated  by a finite set $\hat {\cal A}$ (because it is of finite index in a finitely generated
group $G$). Hence $K$ is generated by the finite set ${\cal C} = \varphi(\hat {\cal A})$.
For every element $x \in {\cal C}$ choose
one element $y \in H$ from its preimage under $\varphi$ and denote by $\bar {\cal C} $ the subset of $H$ consisting of
them. Since
$$H / M \cong K~,$$ $H$ is generated by the finite set ${\cal B}= \bar {\cal C} \cup M$.

Now one can define the corresponding length functions $|\cdot|_H$ and $|\cdot|_{\G}$ which satisfy the following
properties:
$$\forall~g \in {\G} \quad |\varphi(g)|_K \le |g|_{\G}~, $$
$$\forall~h\in H\quad |h|_H \le |\varphi(h)|_K+1~.$$
Combining these inequalities we get
$$\forall~h\in H\quad |h|_H \le |h|_{\G}+1~. $$
Therefore, $D_H(n) \le n+1$, i.e. $H$ is undistorted in $\G$. The group $G$ is hyperbolic and
any its subgroup of finite index is quasiconvex, thus, by lemma \ref{quasi-undist}, $\G$ is undistorted in $G$.
Evidently the property that a subgroup is undistorted in a group is transitive, hence $H$ is undistorted in $G$.
It remains to apply lemma \ref{quasi-undist}
to complete the proof. $\square$

\vspace{.15cm}
\underline{Proof of theorem \ref{ness}.} It is enough to show that $M=ker \phi \cap H$ is finite.
By the conditions,
$$M \subset P_1Q^{-1}QP_2=\bigcup_{x\in P_1,y \in P_2} xQ^{-1}Qy~.$$
Proving by contradiction, assume that $M$ is infinite. Then, since the subsets $P_1,P_2$ are finite,
there are elements $g \in P_1$ and $h \in P_2$ such that the intersection
$$A \stackrel{def}{=}M \cap gQ^{-1}Qh$$ is infinite. Therefore, $B=g^{-1}Ah^{-1}$ is an infinite subset of
$Q^{-1}Q$ satisfying $$\phi(B)=\{\phi(g^{-1}h^{-1})\}~\mbox{-- a one-element subset}~.$$
It is easy to see that the latter is impossible if $\phi$ is a quasiisometry between $Q$ and $\phi(Q)$
(since for any $u,v \in Q$~ $d(u,v)=|u^{-1}v|_G$, $d_1\left(\phi(u),\phi(v)\right)=|\phi(u^{-1}v)|_{G_1}$,
$u^{-1}v \in Q^{-1}Q$). $\square$

\vspace{.15cm}
To prove theorem \ref{qcex} we will need the following statement:
\begin{lemma} \label{finiteindx} {\normalfont \cite[Thm. 2]{paper2}} Assume that $U$ is a finite union of
quasiconvex products in
a hyperbolic group $G$ and the subgroups $F_1,F_2,\dots,F_l$ are all the members of $U$.
If $H$ is a subgroup of $G$ and $H \subseteq U$ then
for some $g \in G$ and $j\in \{1,2,\dots,l\}$ one has $|H:(H\cap gF_jg^{-1})|<\infty$.
\end{lemma}

\vspace{.15cm}
\underline{Proof of theorem \ref{qcex}.} Fix arbitrary finite subsets $P_1,P_2$ of $G$.
Then the set $U=P_1Q^{-1}QP_2$ is also a finite union of quasiconvex products in the group $G$ with the
same members $F_1,\dots,F_l$. Now the claim of the theorem follows directly from lemma \ref{finiteindx}. $\square$

\section{Auxiliary Lemmas}
Assume $\cal X$ is a $\delta$-hyperbolic metric space with metric $d(\cdot,\cdot)$.

\begin{lemma} \label{longpath} Suppose $\varkappa \ge 0$, $X,Y,Z,X',Y' \in {\cal X}$ and  $X' \in {\cal O}_\varkappa([X,Z])$,
$Y' \in {\cal O}_\varkappa([Y,Z])$. Then $(X'|Y')_Z \le (X|Y)_Z +2\varkappa$.
\end{lemma}

\underline{Proof.} Let $X'' \in [X,Z]$, $Y'' \in [Y,Z]$ satisfy $d(X',X'') \le \varkappa$, $d(Y',Y'') \le \varkappa$.
By the triangle inequality
$$(X'|Y')_Z = \frac12 \bigl( d(X',Z)+d(Y',Z)-d(X',Y') \bigr) \le $$
$$\le \frac12 \bigl( d(X'',Z)+d(Y'',Z)-d(X'',Y'')+2d(X',X'')+2d(Y',Y'') \bigr)~. $$

Now, since $d(X,Z)=d(X'',Z)+d(X'',X)$ and $d(Y,Z)=d(Y'',Z)+d(Y'',Y)$, we achieve
$$(X'|Y')_Z \le \frac12 \bigl( d(X,Z)+ d(Y,Z) -[d(X'',Y'')+d(X'',X)+d(Y'',Y)]\bigr) + 2\varkappa \le$$
$$\le \frac12 \bigl( d(X,Z)+ d(Y,Z) - d(X,Y)\bigr) + 2\varkappa = (X|Y)_Z +2\varkappa~.$$
Q.e.d $\square$

\begin{lemma} \label{quasigeod} Let $\bar \lambda >0$, $\bar c \ge 0$, $C_0 \ge 14 \delta$,
 $C_1 = 12(C_0+\delta)+\bar c + 1$ be given.
Then  for $\lambda = \bar \lambda/4>0$ there exist $c=c(\bar \lambda,\bar c,C_0) \ge 0$ satisfying the statement below. \\
Assume $N \in \N$, $X_i \in {\cal X}$, $i=0,\dots,N$, and $q_i$ are $(\bar \lambda,\bar c)$-quasigeodesic paths between
$X_{i-1}$ and $X_i$ in $\cal X$, $i=1,\dots,N$. If $\|q_i\| \ge (C_1+ \bar c)/{\bar \lambda}$, $i=1,\dots,N$, and
$(X_{i-1}|X_{i+1})_{X_i} \le C_0$ for all $i=1,\dots,N-1$,
then the path $q$ obtained as a consequent concatenation of
$q_1,q_2,\dots,q_N$, is $(\lambda,c)$-quasigeodesic.
\end{lemma}

\underline{Proof.} Let the number $\nu=\nu(\bar \lambda,\bar c) \ge 0$ be chosen according to the claim of lemma
\ref{close}. Set $c=\frac52(\nu+C_1)\ge 0$.

Suppose $p$ is an arbitrary subpath of $q$. Then $p_- \in q_j$, $p_+ \in q_k$ for some
$1 \le j \le k \le N$. If $j=k$, $p$ is a subpath of $q_j$ and therefore, it is $(\bar \lambda,\bar c)$-quasigeodesic,
hence it is $(\lambda,c)$-quasigeodesic.

Now let's assume that $j<k$. By our conditions and the choice of $\nu$, there are points $U \in [X_{j-1},X_j]$,
$V \in [X_{k-1},X_k]$ such that $d(p_-,U) \le \nu$ and $d(p_+,V) \le \nu$.
$\|[X_{j-1},X_j]\| \ge \bar \lambda \|q_j\| - \bar c \ge C_1$, similarly, $[X_{k-1},X_k] \ge C_1$, hence after
shifting the points $U$ and $V$ along the segments $[X_{j-1},X_j]$ and $[X_{k-1},X_k]$ (correspondingly) by distances
at most $C_1$ we will obtain $\|[U,X_j]\| \ge C_1$, $\|[X_{k-1},V]\| \ge C_1$, $d(U,p_-) \le \nu + C_1$,
$d(V,p_+)\le \nu + C_1$ and, therefore, $d(p_-,p_+) \ge d(U,V)-2\nu-2C_1$.

According to the lemma \ref{longpath}, all the conditions of lemma \ref{brokenlines2} applied to the broken line
$[U,X_j,\dots,X_{k-1},V]$ are satisfied, hence $$d(U,V) \ge \frac12\|[U,X_j,\dots,X_{k-1},V]\|~.$$
Consequently, $$d(p_-,p_+) \ge \frac12 \left(d(U,X_j)+ \sum_{i=j}^{k-2} d(X_i,X_{i+1})+d(X_{k-1},V) \right)
-2\nu-2C_1~.$$

Finally, we observe that $\displaystyle d(U,X_j) \ge \frac12 d(U,X_j) + \frac{C_1}2>\frac12(d(U,X_j)
+ \bar c)$ and analogously for the
other summands. Denote by $q_j'$, $q_{k}'$ the segments of $q_j$ and $q_{k}$ from $p_-$ to $X_j$ and
from $X_{k-1}$ to $p_+$ correspondingly. We obtain

$$d(p_-,p_+) \ge \frac14 \left(d(U,X_j) +\bar c+ \sum_{i=j}^{k-2} \Bigl(d(X_i,X_{i+1})+\bar c\Bigr)+
d(X_{k-1},V)+\bar c \right) -$$ $$-2\nu-2C_1 \ge$$

$$ \ge \frac14 \left(d(p_-,X_j)+\bar c+ \sum_{i=j}^{k-2} \Bigl(d(X_i,X_{i+1})+\bar c\Bigr)+d(X_{k-1},p_+)
+\bar c\right) -\frac52\nu-\frac52C_1 \ge$$
$$\ge \frac14 \left(\bar \lambda \|q_j'\|+ \sum_{i=j}^{k-2} \bar \lambda \|q_{i+1}\|+ \bar \lambda\|q_k'\|
\right) -\frac52\nu-\frac52C_1 \ge \frac{\bar \lambda}4 \|p\|-\frac52(\nu+C_1)~.$$
The statement is proved. $\square$

Now let $G$ be a $\delta$-hyperbolic group, $\delta \ge 0$.

\begin{lemma} \label{y-mod} Let $H$ be a non-elementary subgroup of a hyperbolic group $G$, and
$g$ be an $H$-suitable element. If $y \in C_H\bigl(E(H)\bigr) \backslash E(g)$ then there exists $N \in \N$
such that the element $yg^n$ has infinite order in $H$ and is $H$-suitable for every $n\ge N$.
\end{lemma}

\underline{Proof.}  Observe that
$E(g) \neq E(ygy^{-1})$ because, otherwise, we would have $yg^ky^{-1}=g^{l}$ for some non-zero integers
$k,l$, and (\ref{elemdef}) would imply $y \in E(g)$ which is not true by the conditions of the lemma.
Therefore $E(g)\cap E(ygy^{-1})$ is finite, hence, by lemma \ref{quasintersection}
there is $C_0' \ge 0$ such that
$(g^l|yg^ky^{-1})_{1_G} \le C_0'$ for any $k,l \in \Z$. Set $C_0=C_0'+2|y|_G+14\delta$, then
$C_0 \ge 14 \delta$ and
\begin{equation} \label{smallprod} (g^{-n}y^{-1}|yg^n)_{1_G} \le (g^{-n}|yg^ny^{-1})_{1_G} +2|y|_G \le C_0 \quad
\forall~n \in \N.
\end{equation}

First, let us show that $yg^n \in H^0$ if $n\in \N$ is sufficiently big.

Choose $C_1=12(C_0+\delta)+1$ and $N_1\in \N$ so that $|yg^n|_G \ge C_1$ for all $n \ge N_1$.
Suppose $(yg^n)^t=1_G$ for some $t \in \N$ and $n \ge N_1$. Consider the broken line $[X_0,X_1,\dots,X_t]$
in \ga with $X_i=(yg^n)^i$, $i=0,1\dots,t$. We can estimate
$$ \bigl(X_{i-1}|X_{i+1}\bigr)_{X_i} =
\bigl((yg^n)^{i-1}|(yg^n)^{i+1}\bigr)_{(yg^n)^{i}}=\bigl(g^{-n}y^{-1}|yg^n\bigr)_{1_G} \le C_0~. $$

It satisfies all the conditions of lemma \ref{brokenlines2}, therefore
$$\|[X_0,X_t]\| \ge \frac12 \|[X_0,X_1,\dots,X_t]\|\ge C_1/2>0~,$$ but we assumed that $X_0=X_t$. A contradiction.
Hence the element $yg^n$ has infinite order for each $n \ge N_1$.

Now, if $n \ge N_1$, $E(yg^n) \neq E(g)$ because, otherwise, we would obtain $yg^n \in E(g)$ which implies
$y \in E(g)$ contradicting to the conditions of the lemma. Hence,
$$E(g) \cap E(yg^n) \subset T(g) = E(H) \subset E(yg^n), ~\mbox{ thus }~ E(g) \cap E(yg^n)=E(H).$$

Let $w_1,w_2$ be shortest words in the alphabet $\cal A$ representing $y$ and $g$ correspondingly.
By lemma \ref{power} there exist $\bar \lambda>0$ and $\bar c' \ge 0$ such that any path in \ga labelled
by the word $w_2^n$ is $(\bar \lambda,\bar c')$-quasigeodesic for any $n \in \N$. Consequently, any path labelled
by $w_1w_2^n$ is $(\bar \lambda,\bar c)$-quasigeodesic where $\bar c =\bar c'+2\|w_1\|$. Set
$C_1 = 12(C_0+\delta)+\bar c + 1$. Suppose that $n \ge (C_1+\bar c)/{\bar \lambda}$.
Then $\|w_1w_2^n\|\ge (C_1+\bar c)/{\bar \lambda}$ and by (\ref{smallprod}) we can apply
lemma \ref{quasigeod} to find $\lambda>0$ and $c\ge 0$ (not depending on $n$) such that any path labelled by
$(w_1w_2^n)^t$ is $(\lambda,c)$-quasigeodesic for any $t \in \Z$.

Suppose $x \in E(yg^n)$. There is $k \in \N$ such that $x(yg^n)^kx^{-1} = (yg^n)^{\epsilon k}$ where
$\epsilon \in \{1,-1\}$, hence $x(yg^n)^{lk}x^{-1} = (yg^n)^{\epsilon lk}$ for any $l \in \N$.
Consider the geodesic quadrangle $Y_1Y_2Y_3Y_4$ in \ga with $Y_1=1_G$, $Y_2=x$, $Y_3=x(yg^n)^{lk}$,
$Y_4=x(yg^n)^{lk}x^{-1}$ and $(\lambda,c)$-quasigeodesic paths $p$ between $Y_2$ and $Y_3$
and $q$ between $Y_1$ and $Y_4$ labelled by the words $(w_1w_2^n)^{lk}$ and $(w_1w_2^n)^{\epsilon lk}$
correspondingly. Choose $\nu=\nu(\lambda,c)$ to be the constant given by lemma \ref{close}. Thus
$$p \subset {\cal O}_{\nu}([Y_2,Y_3]),~[Y_2,Y_3] \subset {\cal O}_{\nu}(p),~
q \subset {\cal O}_{\nu}([Y_1,Y_4]),~[Y_1,Y_4] \subset {\cal O}_{\nu}(q)~.$$

Obviously, by taking the number $l$ sufficiently large, one can find a subpath $r$ of $p$ labelled by
$w_2^n$ with its endpoints $r_-$ and $r_+$ having distances at least $(|x|_G+\nu)$ from both of the vertices
$Y_2$ and $Y_3$. Then an application of lemma \ref{quadrangle} will give us
$$r_-,r_+ \in {\cal O}_{\nu+2\delta}([Y_1,Y_4]) \subset {\cal O}_{2\nu+2\delta}(q)~.$$
Let denote $u,v$ denote the points on the path $q$ with $d(r_-,u) \le 2\nu+2\delta$ and $d(r_+,v) \le 2\nu+2\delta$,
and let $r'$ be the subpath of $q$ (or $q^{-1}$) starting at $u$, ending at $v$. The length of $r$
(and, hence, the length of $[r_-,r_+]$) depends on $n$, thus it can be made as large as we please,
therefore $r'$ will also be long (compared to $\|w_1\|$), consequently $r'$ will have a subpath $q'$ labelled
by $w_2^t$, $t \in \Z$, with $|t| \ge \bar \lambda n/3$. Since the quadrangles in \ga are $2\delta$-slim, we achieve
$$[u,v]\subset {\cal O}_{2\delta}([r_-,r_+]\cup [r_-,u] \cup [r_+,v])\subset {\cal O}_{2\nu+4\delta}([r_-,r_+])~,
~\mbox{ hence}$$
$$q' \subset {\cal O}_{\nu}([u,v]) \subset {\cal O}_{3\nu+4\delta}([r_-,r_+]) \subset {\cal O}_{4\nu+4\delta}(r)~.$$

Consider the vertices $a_0=q'_-,a_2,\dots,a_{|t|}=q'_+$ of the path $q'$ such that the subpaths between $a_{i-1}$
and $a_i$ are labelled by $w_2$ (respectively, $w_2^{-1}$ if $t<0$) for every $1\le i \le |t|$
(we will call them {\it phase} vertices). Then each
of them is at distance at most $(4\nu+4\delta+\|w_2\|)$ from some phase vertex of $r$. There are only finitely many
words over the alphabet $\cal A$ of length at most $(4\nu+4\delta+\|w_2\|)$, therefore, if $n$ is sufficiently
large, there will be two paths $\alpha$ and $\beta$ connecting two different phase vertices of $q'$ with
some vertices of $r$ having the same word $w_3$ written on them. Thus we achieve an equality in the group $G$:
$$w_2^s=w_3w_2^{s'}w_3^{-1} ~\mbox{ for some } s,s' \in \Z\backslash \{0\}~.$$

So, if $z$ denotes the element of $G$ represented by the word $w_3$, we have
\begin{equation} \label{z} zg^{s'}z^{-1} = g^s
\end{equation}
Then $z \in E(g)$ by (1). By the construction,
$x=(yg^n)^{\zeta}g^\xi z g^{\xi'}(yg^n)^{\zeta'}$ for some $\zeta,\xi,\zeta',\xi' \in \Z$. Note that
$(yg^n)^{-\zeta} x (yg^n)^{-\zeta'} \in E(yg^n)$ and $g^\xi z g^{\xi'} \in E(g)$. Hence,
$$g^\xi z g^{\xi'}= (yg^n)^{-\zeta} x (yg^n)^{-\zeta'} \in E(yg^n) \cap E(g) =E(H)~.$$
Now, since $yg^n \in C_H\bigl( E(H)\bigr)$ we obtain  $x \in \langle yg^n \rangle \cdot E(H) =
 \langle yg^n \rangle \times E(H)$ for arbitrary $x$ from $E(yg^n)$.
 This implies that $yg^n$ is $H$-suitable. $\square$

\begin{lemma} \label{y-modnon-comm} Let $G$ be a hyperbolic group, $k \in \N$ and let
$g_1,g_2,\dots,g_k $ be pairwise non-commensurable elements of $G$. 
Consider $y_i \in G \backslash E(g_i)$ for each $i=1,2,\dots,k$.
Then there exists $N \in \N$
such that the elements $y_1g_1^n,\dots,y_kg_k^n$ have infinite order and are pairwise non-commensurable
if $n\ge N$. 
\end{lemma}

\underline{Proof.} 
By the same argument as above, $y_ig_i^n \in G^0$ for each $i\in \{1,2,\dots,k\}$.
Suppose that $y_ig_i^n$ is commensurable with
$y_jg_j^n$, $1 \le i < j \le k$. Then there is $x \in G$ satisfying $x(y_ig_i^n)^tx^{-1}=(y_jg_j^n)^{t'}$ for some
$t,t' \in \Z\backslash \{0\}$. Therefore,
if $n$ is sufficiently large, we can prove that there exist $z\in G$ and
$s,s' \in \Z\backslash \{0\}$ such that $zg_i^{s'}z^{-1} = g_j^s$
in exactly the same way as we proved (\ref{z}). Which leads to a contradiction with the assumptions of the lemma.
$\square$

\begin{lemma} \label{findK} Suppose $G$ is a hyperbolic group, $H$ is its non-elementary subgroup and
$\alpha_1,\dots,\alpha_n$
are points on the boundary $\partial G$. Then $H$ has a non-elementary subgroup $K$ that is quasiconvex in $G$
and $\alpha_i \notin \Lambda(K) \subset \partial G$ for every \\$i=1,2,\dots,n$.
\end{lemma}

\underline{Proof.} Induction on $n$.

Let $n=1$. Choose elements of infinite order $g_1,g_2 \in H$ with $E(g_1) \neq E(g_2)$ and
let the words $w_1,w_2$ over the alphabet $\cal A$ represent them. Following the proof
of \cite[Cor. 6]{Olsh2}, we get $M \in \N$, $\lambda>0$ and $c \ge 0$ such that any word of the form
$$w_{i_1}^{mk_1}\cdot \dots\cdot w_{i_s}^{mk_s}$$ where $i_j \in \{1,2\}$, $i_j\neq i_{j+1}$,
$k_j \in \Z \backslash \{0\}$ and $m \ge M$, is $(\lambda,c)$-quasigeodesic. Take $\nu=\nu(\lambda,c)$
from lemma \ref{close}.
Then the subgroup $K_1=\langle g_1^m,g_2^m\rangle \le H$ is free of rank $2$ and $\varepsilon$-quasiconvex
in $G$ ($\varepsilon=\nu+m \cdot \max\{\|w_1\|,\|w_2\|\}$). By taking $m$ large enough, we can obtain
$|H:K_1|=\infty$.

Now, if $\alpha_1 \notin \Lambda(K_1)$, there is nothing to prove. So, assume that $\alpha_1 \in \Lambda(K_1)$.
By lemmas \ref{conj} and \ref{newdouble} there exists $h\in H$ such that $card(K_1 \cap hK_1h^{-1})<\infty$.
The subgroup $hK_1h^{-1} \le H$ is non-elementary and quasiconvex in $G$, hence,
using lemma \ref{quasintersection} and the definition of the boundary $\partial G$ we obtain
$$\Lambda(K_1) \cap \Lambda(hK_1h^{-1}) = \emptyset ~\mbox{ in } \partial G~.$$
Consequently, $\alpha_1 \notin \Lambda(hK_1h^{-1})$.

Assume, now, that $n>1$. And the induction hypothesis is verified for $\alpha_1,\dots,\alpha_{n-1} \in \partial G$.
I.e. there is a non-elementary subgroup $K' \le H$ with $\alpha_i \notin \Lambda(K')$, $1\le i \le n-1$.
Using the base of our induction, we obtain a non-elementary subgroup $K \le K' \le H$ that is
quasiconvex in $G$ and $\alpha_n \notin \Lambda(K)$. Since $\Lambda(K) \subseteq \Lambda(K')$, $K$ satisfies
all the conditions needed.

The proof of the lemma is complete. $\square$

It is a well-known fact that the set of all rational points $\{g^\infty~|~g \in G^0\}$ is dense in the group
boundary $\partial G$ (see, for example, \cite[Theorem]{Ruane}, \cite[4,8.2D]{Gromov}).
We will need a bit stronger statement:

\begin{lemma} \label{density} Assume $H$ is a non-elementary subgroup of a hyperbolic group $G$ and
$\alpha \in \partial G$.
Then ~$\Lambda(H) \subseteq cl(H \circ \alpha)$ where $H \circ \alpha$ is the orbit of $\alpha$ under the action of $H$
and $cl(H \circ \alpha)$ is its closure inside of $\partial G$.
\end{lemma}

\underline{Proof.} Since $H$ is non-elementary, the set $H \circ \alpha$ consists of more than one point.
By definition, $H \subset St_G(H \circ \alpha)$, hence, applying lemma \ref{stabilizer}, we achieve
$$\Lambda(H) \subseteq cl(H \circ \alpha)~.$$ Q.e.d. $\square$

\section{On condition $(*)$} \label{*-cond}

Let $G$ be a $\delta$-hyperbolic group, $Q \subseteq G$ --
$\eta$-quasiconvex subset.

\begin{lemma} \label{q-1q}The subset $Q^{-1}Q \subseteq G$ is $(\eta+\delta)$-quasiconvex.
\end{lemma}

\underline{Proof.} Consider arbitrary $x \in Q^{-1}Q$, $x=u^{-1}v $ where $u,v \in Q$. Then
$[u,v] \subset {\cal O}_{\eta}(Q)$. Since the metric on \ga is invariant under the action of
$G$ by left translations, we have
\begin{equation} \label{xuv}
[1_G,x]=u^{-1} \circ [u,v] \subset {\cal O}_{\eta}(u^{-1}Q) \subset {\cal O}_{\eta}(Q^{-1}Q)~.
\end{equation}

Since the geodesic triangles in \ga are $\delta$-slim, for any two $x_1,x_2 \in Q^{-1}Q$ using (\ref{xuv})
one obtains $$[x_1,x_2] \subset {\cal O}_{\delta}([1_G,x_1]\cup[1_G,x_2])
\subseteq {\cal O}_{\delta+\eta}(Q^{-1}Q)~.$$ The lemma is proved. $\square$

\begin{lemma} \label{SQ} Suppose $S,Q \subseteq G$ and the subset $Q$ is $\eta$-quasiconvex. Then on the boundary
$\partial G$ of the group $G$ one has
$$ \Lambda(S \cdot Q) \subseteq \Lambda(S)\, \cup \,(S\cdot Q) \circ \Lambda(Q^{-1} \cdot Q)~,$$
$$ \Lambda(S \cdot Q^{-1}) \subseteq \Lambda(S)\, \cup \,(S\cdot Q^{-1} \cdot Q) \circ \Lambda(Q^{-1})~.$$
\end{lemma}

\underline{Proof.} Let $P \subseteq G$, consider an arbitrary limit point $\alpha \in \Lambda(SP)$. There is a sequence
$(z_i)_{i\in \N}$ converging to infinity in $G$ with $z_i =x_iy_i$, $x_i \in S$, $y_i \in P$ for all $i \in \N$.

%

{\it I.} Suppose, first, that $ \sup_{i \in \N}(z_i|x_i)_{1_G} =\infty$. Then one can find a sequence
$(i_j)_{j\in \N}$ of natural numbers such that
$$\lim_{j \to \infty} (z_{i_j}|x_{i_j})_{1_G} = \infty~.$$
But $\displaystyle \lim_{j \to \infty}z_{i_j} =\lim_{i \to \infty}z_{i}= \alpha$,
which implies that $(x_{i_j})_{j\in \N}$ also converges to
infinity and $$\lim_{j \to \infty}x_{i_j} = \lim_{j \to \infty}z_{i_j} = \alpha~. $$
Thus, $\alpha \in \Lambda(S)$.

{\it II.} Therefore, we can now assume that there is a number $M \ge 0$ such that $(z_i|x_i)_{1_G} \le M$ for every
$i\in \N$.
For each $i \in \N$ consider a geodesic triangle in \ga with vertices $1_G$, $x_i$ and $z_i$. It is $\delta$-thin,
hence $d(1_G,[x_i,z_i])\le M+\delta$.

a). Suppose $P=Q$. Fix an arbitrary element $q \in Q$ and let $\varkappa = |q|_G$. Then
$$[1_G,y_i] \in \O_{\delta}([1_G,q] \cup [q,y_i]) \subset \O_{\delta+\varkappa}([q,y_i])
\subset \O_{\delta+\varkappa+\eta}(Q)~.$$

Using the left translation-invariance of the word metric, we get
$$[x_i,z_i]=[x_i,x_iy_i]\subset \O_{\delta+\varkappa+\eta}(x_iQ)~.$$
Consequently, there exists $q_i \in Q$ satisfying
$$d(1_G,x_iq_i) =|x_iq_i|_G\le M+2\delta+\varkappa+\eta~\mbox{ for every } i \in \N.$$

The group $G$ has only finitely many elements in a ball of finite radius, hence, by passing to a subsequence,
we can assume that $x_iq_i=p \in SQ$ for all $i \in \N$. Thus, $z_i=x_iq_iq_i^{-1}y_i=pq_i^{-1}y_i \in pQ^{-1}Q$
for every $i$, which implies $$\alpha \in \Lambda(pQ^{-1}Q)=p\circ \Lambda(Q^{-1}Q) \subset (SQ)\circ \Lambda(Q^{-1}Q)~.$$

b). Assume, $P=Q^{-1}$. Then $y_i^{-1} \in Q$ hence
$$[x_i,z_i]=[z_iy_i^{-1}, z_i] \subset \O_{\delta+\varkappa+\eta}(z_iQ)~.$$
So, there are elements $q_i \in Q$ such that
$$d(1_G,z_iq_i) =|z_iq_i|_G\le M+2\delta+\varkappa+\eta~\mbox{ for every } i \in \N.$$

As before, we can suppose that $z_iq_i=p \in SQ^{-1}Q$ for all $i \in \N$. Thus $z_i=pq_i^{-1} \in pQ^{-1}$
for every $i$, implying $$\alpha \in \Lambda(pQ^{-1})=p\circ \Lambda(Q^{-1}) \subset (SQ^{-1}Q)\circ \Lambda(Q^{-1})~.$$
Q.e.d. $\square$

\begin{lemma} \label{findalpha} Assume that $H$ is a subgroup  and $A$ is a non-empty quasiconvex
subset of a hyperbolic group $G$. The following conditions are equivalent:

{\bf 1.} There are finite subsets ~$P_1,P_2 \subset G$~  such that~
$ \displaystyle H \subseteq P_1 \cdot A \cdot P_2~;$

{\bf 2.} $ \displaystyle  \Lambda(H) \subseteq G \circ \Lambda(A).$

\end{lemma}

\underline{Proof.}  The implication 1 $\Rightarrow$ 2 is an immediate consequence of lemma \ref{limitsets}
since $$\Lambda(P_1 \cdot A \cdot P_2)=\Lambda(P_1A)=P_1\circ \Lambda(A) \subset G \circ \Lambda(A)~.$$
Now let's show that 2 implies 1. Denote $\Omega=\Lambda(A)$.

If the subgroup $H$ is finite then the claim is trivial.

If $H$ is infinite elementary then $card\bigl(\Lambda(H)\bigr)= 2$, hence, according to the condition $2$, there are
elements $g_1,g_2 \in G$ such that
$$\Lambda(H) \subset g_1 \circ \Lambda(A) \cup g_2 \circ \Lambda(A)= \Lambda(g_1A \cup g_2A)~.$$
The subset $g_1A\cup g_2A \subset G$ is quasiconvex by lemma \ref{q-cunion,product}, hence we can apply lemma
\ref{tame-qc} to get a finite subset $P_2$ of $G$ satisfying
$H \subseteq (g_1A\cup g_2A)P_2=P_1AP_2$ where $P_1=\{g_1,g_2\}$.

Thus we can assume that the subgroup $H$ is non-elementary.

\noindent {\it Case 1.}
Suppose that for some $g \in G$~
$g\circ\Omega$ contains a non-empty open set $U$ of the subspace $\Lambda(H)$, i.e. $U=U'\cap \Lambda(H)$
for some open set $U' \subseteq \partial G$. Then, by lemma \ref{density},
for any $\beta \in \Lambda(H)$ there exists $h \in H$ such that $h \circ \beta \in U'$. On the other hand,
$h \circ \beta \in \Lambda(H)$ by lemma \ref{invariant}, thus, $h \circ \beta \in U$, i.e.
$\beta \in h^{-1} \circ U$. Consequently,
\begin{equation}  \label{cover} \Lambda(H) \subseteq \bigcup_{h \in H} h \circ U~.
\end{equation}

The space $\Lambda(H)$ is a closed subspace of the compact metric space $\partial G$, hence it is compact
itself and one can choose a finite subcover of the open cover from (\ref{cover}). Thus
$$\Lambda(H) \subseteq \bigcup_{i=1}^N h_i \circ U \subseteq \bigcup_{i=1}^N h_i \circ (g\circ\Omega)
= \bigcup_{i=1}^N h_i \circ \Lambda \bigl(gA\bigr) = \Lambda\left(\bigcup_{i=1}^N h_igA\right)=\Lambda (P_1A) $$
according to lemma \ref{limitsets}, where $\displaystyle P_1=\bigcup_{i=1}^N h_ig \subset G$, $card(P_1)<\infty$.

The set ~$\displaystyle P_1 A=\bigcup_{y\in P_1} yA$~ is quasiconvex as a finite union of quasiconvex sets,
therefore, we can apply
lemma \ref{tame-qc} to find a finite subset $P_2$ of the group $G$ such that
$H \subseteq P_1 \cdot A \cdot P_2$ as we needed.

Hence, we can proceed to \\
{\it Case 2.} For every $g \in G$ ~$g\circ \Omega$ contains no non-empty open subsets of $\Lambda(H)$.
$\Omega$ is a closed subset of the boundary $\partial G$ by lemma \ref{limitsets}.(b), thus
$g \circ \Omega$ is also closed and, hence, $\bigl(g\circ \Omega \bigr)\cap \Lambda(H)$ is a closed nowhere dense subset
of the compact metric space $\Lambda(H)$. Evidently, $\Lambda(H)$ is a Baire space (it is locally compact and Hausdorff).
Since the group $G$ is countable, the set
$$(G\circ \Omega) \cap \Lambda(H) = \bigcup_{g \in G} (g \circ \Omega) \cap \Lambda(H)$$ is of the first category
in the space $\Lambda(H)$, hence, by a well-know theorem from
topology (see, for instance, \cite[Chap. $XI$, Thm. 10.5]{Dugundji}),
$$\Lambda(H) \neq (G\circ \Omega) \cap \Lambda(H)~,$$
 therefore, $\Lambda(H) \not \subseteq G \circ \Lambda(A)$ which is a contradiction
 to our assumptions. Thus, Case $2$ is impossible. $\square$

\vspace{.15cm}
\remark \label{nowheredense} Assume (in the notations of lemma \ref{findalpha}) that $H$ is non-elementary.
Then the following two properties are equivalent:

{\bf 1.} There are no finite subsets $P_1,P_2$ of $G$ such  that $H \subset P_1AP_2$;

{\bf 2.} On the hyperbolic boundary $\partial G$ for any $g \in G$ the set
$\bigl(g \circ\Lambda(A)\bigr) \cap \Lambda(H)$ is nowhere dense in $\Lambda(H)$.

\vspace{.15cm}
In the proof of lemma \ref{findalpha} the condition 1 automatically puts us into the Case 2, thus 1~$\Rightarrow$~2.
Now, if the property 2 holds and the property 1 doesn't, we can find some finite subsets $P_1,P_2 \subset G$ satisfying
$H \subset P_1AP_2$. Therefore, by lemma \ref{limitsets},
$$\Lambda(H)\subset\Lambda(P_1AP_2)=\Lambda(P_1A)~.$$
Hence $\displaystyle \Lambda(H)=\bigcup_{g \in P_1} \bigl(g \circ \Lambda (A) \cap \Lambda(H)\bigr)$ contradicting to
2 because of the fact that a finite union of nowhere dense subsets is nowhere dense in $\Lambda(H)$. Hence
2~$\Rightarrow$~1.

\begin{lemma} \label{smallunion} Suppose that $H$ is a non-elementary subgroup and $Q,S$ are quasiconvex subsets
of a hyperbolic group $G$. Assume that for any two finite subsets $P_1,P_2$ of the group $G$
\begin{equation} \label{Q,S} H \nsubseteq P_1Q^{-1}QP_2~ \mbox{and}~H \nsubseteq P_1S^{-1}SP_2~. \end{equation}
Then the (quasiconvex) subsets $T_1=Q\cup S$ and $T_2=QS$ satisfy the same property: for any $i \in \{1,2\}$ and
arbitrary finite $P_1,P_2 \subset G$ one has $$H \nsubseteq P_1T_i^{-1}T_iP_2~.$$
\end{lemma}

\underline{Proof.} a). Since $T_1^{-1}=Q^{-1} \cup S^{-1}$, we can use lemmas \ref{limitsets}.(c) and \ref{SQ} to obtain
$$\Lambda\left(T_1^{-1}T_1\right)=\Lambda\left(Q^{-1}Q \cup Q^{-1}S \cup S^{-1}Q \cup S^{-1}S\right)= $$
$$=\Lambda\left(Q^{-1}Q \right) \cup \Lambda \left( Q^{-1}S \right) \cup \Lambda \left(S^{-1}Q \right) \cup
\Lambda \left(S^{-1}S\right)\subseteq$$
$$\subseteq \Lambda\left(Q^{-1}Q \right) \cup\Lambda(Q^{-1}) \cup G\circ\Lambda \left( S^{-1}S \right) \cup
\Lambda(S^{-1}) \cup G \circ \Lambda \left(Q^{-1}Q \right) \cup \Lambda \left(S^{-1}S\right) =$$
$$=G\circ\Lambda \left( S^{-1}S \right) \cup G\circ\Lambda \left( Q^{-1}Q \right)=
G\circ\Lambda \left( S^{-1}S \cup Q^{-1}Q \right)$$
(here we used the fact that if $s \in S$ then $S^{-1}s \subset S^{-1}S$, and by lemma \ref{limitsets}.(d)
$\displaystyle \Lambda(S^{-1}) = \Lambda(S^{-1}s) \subset \Lambda(S^{-1}S)$ and similarly for $Q$~).\\
Thus, $\displaystyle G\circ\Lambda\left(T_1^{-1}T_1\right) \subseteq G\circ\Lambda \left( S^{-1}S \cup Q^{-1}Q \right)$.

The conditions (\ref{Q,S}) imply (by remark \ref{nowheredense}) that for any $g \in G$ the sets
\\ $g \circ\Lambda(Q^{-1}Q)\cap \Lambda(H)$ and $g \circ\Lambda(S^{-1}S)\cap \Lambda(H)$ are nowhere dense in
$\Lambda(H)$, therefore the set $G \circ \Lambda(Q^{-1}Q \cup S^{-1}S) \cap \Lambda(H)$ is of the first category
in the compact metric space $\Lambda(H)$. Consequently,
\begin{equation} \label{qqss} \Lambda(H) \nsubseteq G \circ (Q^{-1}Q \cup S^{-1}S). \end{equation}
Hence $$\Lambda(H) \nsubseteq G \circ (T_1^{-1}T_1)~.$$ The subset $T_1 \subset G$ is quasiconvex
according to lemma \ref{q-cunion,product}.
Therefore, to finish the proof it remains to apply lemma \ref{q-1q} to $T_1^{-1}T_1$ and then
lemma \ref{findalpha} to $T_1^{-1}T_1$ and $H$.

b). The proof for $T_2$ is similar. Note that $T_2^{-1}=S^{-1}Q^{-1}$, hence, by lemma \ref{SQ}
$$\Lambda \left(T_2^{-1}T_2\right)=\Lambda\left(S^{-1}Q^{-1}QS\right)\subset
\Lambda\left(S^{-1}Q^{-1}Q\right) \cup G\circ \Lambda \left( S^{-1}S \right) ~.$$
$Q^{-1}Q \subset G$ is quasiconvex by lemma \ref{q-1q} and $(Q^{-1}Q)^{-1}= Q^{-1}Q$ therefore, applying lemma
\ref{SQ} two more times we obtain
$$\Lambda \left(T_2^{-1}T_2\right)\subset
\Lambda\left(S^{-1}\right) \cup G\circ\Lambda\left(Q^{-1}Q\right) \cup G\circ \Lambda \left( S^{-1}S \right)=
G \circ \left(Q^{-1}Q \cup S^{-1}S\right) ~.$$
Consequently, recalling (\ref{qqss}), one gets $\Lambda(H) \nsubseteq G \circ (T_2^{-1}T_2)$. Since $T_2$
is a quasiconvex subset of $G$ (lemma \ref{q-cunion,product}), $T_2$ and $H$ satisfy the needed property
by lemma \ref{findalpha}.
$\square$

\begin{cor} \label{q-1q=q} Let $Q$ be a quasiconvex subset of a hyperbolic group $G$ and $H$ be a non-elementary
subgroup of $G$. Assume, in addition, that $Q^{-1} \subset G$ is also quasiconvex. Then the following properties
are equivalent:

{\bf 1.} For arbitrary finite subsets $P_1$, $P_2$ of $G$, $H \nsubseteq P_1QP_2$;

{\bf 2.} For arbitrary finite subsets $P_1$, $P_2$ of $G$, $H \nsubseteq P_1Q^{-1}QP_2$.
\end{cor}

\underline{Proof.} Evidently, $2$ implies $1$. So, let's assume that $1$ holds and prove $2$.
Since the subset $Q^{-1}$ is quasiconvex, we are able to apply lemma \ref{SQ} to achieve
$$\Lambda\left(Q^{-1}Q\right) \subset \Lambda\left(Q^{-1}\right) \cup G\circ\Lambda\left(Q\right)~.$$

Thus, $G \circ \Lambda\left(Q^{-1}Q\right) \subset G\circ\left(\Lambda(Q^{-1}) \cup \Lambda(Q)\right)$.
Observe that the property $1$ is equivalent to $H \nsubseteq P_1Q^{-1}P_2$ for any finite
$P_1,P_2 \subset G$ (because $H^{-1}=H$). Consequently, by remark \ref{nowheredense},
$$\Lambda(H) \nsubseteq G\circ\left(\Lambda(Q^{-1}) \cup \Lambda(Q)\right)~,~ \mbox{ hence ~}
\Lambda(H) \nsubseteq G\circ\Lambda\left(Q^{-1}Q\right)~.$$
Now, by applying lemma \ref{SQ}, we can conclude that the property $2$ holds. $\square$

\underline{\bf Example.} We observe that the implication $1 \Rightarrow 2$ in the latter corollary may fail if
$Q^{-1}$ is not quasiconvex: let $G=F(x,y)$ be the free group with free generators $x,y$. Set $Q$ to be the set of all
reduced words $w$ over the alphabet $\{x^{\pm 1},y^{\pm 1}\}$ satisfying the property:
if $k \in \N$ and $2^k \le \|w\|$ then the letter on $2^k$-th place in $w$ is $x$. Thus,
$$Q=\{x,x^{-1},y,y^{-1},x^2,yx,y^{-1}x,x^3,x^2y,x^2y^{-1},yxy,$$
$$yxy^{-1},y^{-1}xy,y^{-1}xy^{-1},\dots \} \subset G~.$$
The subset $Q$ is quasiconvex in $G$ since any prefix of a word from $Q$ belongs to $Q$.
It is not difficult to show that $G \nsubseteq P_1QP_2$ for any finite subsets $P_1,P_2 \subset G$; nevertheless,
$G=Q^{-1}Q$ (because any reduced word is a suffix of some word from $Q$).

\section{Small Cancellations over Hyperbolic Groups}
In this section we list the results from \cite{Olsh2} that provide the
main tool for proving theorem \ref{mainthm}. 

Let $\cal R$ be a {\it symmetrized} set of words in an alphabet $\cal A$, i.e. if $R \in {\cal R}$ then
$R^{-1}$ and any cyclic permutation of $R$ belong to $\cal R$.

Suppose $G$ is a group with generating set $\cal A$ and $\cal R$ is a symmetrized set of words over $\cal A$.
For the given constants $\varepsilon \ge 0$, $\mu >0$, $0<\lambda \le 1$, $c \ge 0$, $\rho >0$ ~one defines
the {\it generalized  small cancellation conditions} $C(\varepsilon,\mu,\lambda,c,\rho)$,
$C_1(\varepsilon,\mu,\lambda,c,\rho)$, $C_2(\varepsilon,\mu,\lambda,c,\rho)$ and $C_3(\varepsilon,\mu,\lambda,c,\rho)$
(see \cite[Chapter 4]{Olsh2}).

\vspace{.15cm}
\remark \label{indepoflambdac} In the definition of the generalized small cancellation conditions
$C(\varepsilon,\mu,\lambda,c,\rho)$ and $C_j(\varepsilon,\mu,\lambda,c,\rho)$, $j=1,2,3$,
$\lambda,c$ appear only in the condition on every word $R \in {\cal R}$ to be $(\lambda,c)$-quasigeodesic
(a word $R$ is called $(\lambda,c)$-quasigeodesic if any path in \ga labelled by $R$ is $(\lambda,c)$-
quasigeodesic). The number $\rho$ stands for a minimal length of a relation from $\cal R$; $\mu$ has an analogous
meaning to the constant in the classical small cancellation condition $C'(\mu)$;  subwords
of defining relations are distinguished up to factors of length at most $\varepsilon$.

\vspace{.15cm}
Consider non-elementary subgroups $H_1,\dots,H_k$ (some of them may coincide) of a hyperbolic group $G$,
elements $g_i \in H_i$ chosen according to the claim of lemma \ref{lemma3.8} and elements
$x_{i0} \in C_G(E(H_i)) \backslash E(g_i)$, $i=1,\dots,k$. Let $g_i$, $x_{i0},\dots,x_{il}$ be represented by
words $W_i$, $X_{i0},\dots,X_{il}$ over the alphabet $\cal A$ of minimal length, $i=1,\dots,k$.

As the system
${\cal R}={\cal R}_{k,l,m}(W_1,\dots,W_k,X_{10},\dots,X_{kl},m)$ consider all cyclic permutations of $R_i^{\pm 1}$
where
$$R_i \equiv X_{i0}W_i^mX_{i1}W_i^m \dots X_{il}W_i^m~,~~ i=1,2,\dots,k~.$$

Then we have

\begin{lemma} \label{lemma4.2} {\normalfont (\cite[Lemma 4.2]{Olsh2})}
For the words $W_1,\dots,W_k$, $X_{10},\dots,X_{k0}$ given above, there exist $\lambda >0$ such that for any $\mu >0$
there are $l \in \N$ and $c \ge 0$ such that for any $\varepsilon \ge 0$, $\rho>0$ there are $m_0 \in \N$ and words
$X_{11},\dots,X_{kl}$ such that the system ${\cal R}_{k,l,m}$ satisfies $C(\varepsilon,\mu,\lambda,c,\rho)$ and
$C_1(\varepsilon,\mu,\lambda,c,\rho)$-conditions if $m \ge m_0$.
\end{lemma}

\remark \label{remafter4.2} Moreover, from the proof of this lemma it follows that  the elements
$x_{i1},\dots,x_{il} \in H_i$ can be
chosen right after the choice of $l \in \N$ to be any elements that satisfy properties 1)-4) of lemma
\ref{lemma3.7} for $g=g_i$ and $H=H_i$, $i=1,2,\dots,k$.

\vspace{.15cm}

In this article we assume that the concepts of a {\it Van Kampen} ({\it circular}) {\it diagram}
and a {\it Schupp} ({\it annular}) {\it diagram} over a group presentation are known to the reader
(see, for instance, \cite{L-S}).
Let $\cal O$ denote the system of all relations (not only defining) in the group $G$. Let $\cal R$ be some symmetrized
set of additional relations over the alphabet $\cal A$. The group $G_1$ will be defined by its presentation:
\begin{equation} \label{G_1} G_1=\langle {\cal A}~\|~{\cal O} \cup {\cal R} \rangle~. \end{equation}
Thus $G_1$ is a quotient of the group $G$ by the subgroup ${\cal N}=\langle {\cal R}^G\rangle$ that is
a normal closure of the set of elements in $G$ represented by words from $\cal R$.

Inheriting the terminology from \cite{Olsh2} the faces of a Van Kampen diagram with boundary labels
from $\cal O$ (from $\cal R$) will be called $0$-faces ($\cal R$-faces).

In \cite[Chapter 5]{Olsh2} A. Ol'shanskii also introduces {\it reduced} diagrams
(the number of $\cal R$-faces in them can not be reduced after a finite number of certain elementary transformations;
in particular, diagrams with minimal number of $\cal R$-faces are reduced). Further we will only need to know that
for any word $W$ that is trivial in the group $G_1$ there exists a reduced circular diagram over the presentation
(\ref{G_1}) whose boundary label is letter-to-letter equal to $W$. And if two words $U$ and $V$ are conjugate
in $G_1$ then there is a reduced annular diagram over the presentation (\ref{G_1}) whose boundary
contours have labels (letter-to-letter) equal to $U$ and $V$ respectively.

Later in this paper we will consider diagrams over $G$ and $G_1$ (with the presentations
$G=\langle {\cal A}~\|~{\cal O}\rangle$ and (\ref{G_1})), so all of them will be over the alphabet $\cal A$.
A path $q$ inside such a diagram will be called $(\lambda,c)$-{\it quasigeodesic} if a (any) path in the Cayley graph
\ga of the group $G$ with the same label as $q$ is $(\lambda,c)$-quasigeodesic.

The boundary $\partial \Delta$ of a diagram will be divided into at most $4$ distinguished
subpaths (called {\it sections}) each of which will be $(\lambda,c)-$quasigeodesic (for some given $\lambda>0$,
$c \ge 0$).

Suppose $\varepsilon \ge 0$ is a given number.
Consider a simple closed path $o=p_1q_1p_2q_2$ in a diagram $\Delta$ over $G_1$, such that $q_1$ is a subpath of
the boundary cycle of an $\cal R$-cell $\Pi$ and $q_2$ is a subpath of a section $q$ of $\partial \Delta$.
Let $\Gamma$ denote the subdiagram of $\Delta$ bounded by $o$. Assuming that $\Gamma$ has no holes, no
$\cal R$-faces and $\|p_1\|,\|p_2\| \le \varepsilon$, it will be called an $\varepsilon$-{\it contiguity subdiagram}
of $\Pi$ to
$q$. The ratio $\|q_1\|/\|\partial \Pi\|$ will be called the {\it contiguity degree} of $\Pi$ to $q$ and denoted
$(\Pi,\Gamma,q)$.

The following analog of Grindlinger's lemma is proved in \cite[Lemma 6.6]{Olsh2}
(here we include a correction mentioned in \cite{Olsh-SQ}):

\begin{lemma} \label{lemma6.6} For any hyperbolic group $G$ and any $\lambda >0$ there is $\mu_0 >0$ such that
for any $\mu \in (0,\mu_0]$ and $c \ge 0$ there are $\varepsilon \ge 0$ and $\rho >0$ with the following property:

Let the symmetrized presentation (\ref{G_1}) satisfy the $C(\varepsilon,\mu,\lambda,c,\rho)$-condition. Furthermore,
let $\Delta$ be a reduced circular diagram over $G_1$ whose boundary is decomposed into a product of
$(\lambda,c)$-quasigeodesic sections $q^1,\dots,q^r$ where $1\le r \le 4$. Then, provided $\Delta$ has an $\cal R$-face,
there exists an $\cal R$-face $\Pi$ in $\Delta$ and disjoint $\varepsilon$-contiguity subdiagrams
$\Gamma_1,\dots,\Gamma_r$ (some of them may be absent) of $\Pi$ to $q^1,\dots,q^r$ respectively, such that
\begin{equation}  \label{contigsum} (\Pi,\Gamma_1,q^1)+\dots+(\Pi,\Gamma_r,q^r) > 1-23\mu~.
\end{equation}
\end{lemma}

The next lemma is an analog of the previous one for annular diagrams.

\begin{lemma} {\normalfont \cite[Lemma 8.1]{Olsh2}} \label{lemma8.1} For any hyperbolic group $G$ and any $\lambda >0$ there is $\mu_0 >0$ such that
for any $\mu \in (0,\mu_0]$ and $c \ge 0$ there are $\varepsilon \ge 0$ and $\rho >0$ with the following property:

Let the symmetrized presentation (\ref{G_1}) satisfy the $C_1(\varepsilon,\mu,\lambda,c,\rho)$-condition. Further,
let $\Delta$ be a reduced annular diagram over $G_1$ with boundary contours $p=p_1p_2$, $q=q_1q_2$ such that
$p_1$, $p_2$, $q_1$, $q_2$ are $(\lambda,c)$-quasigeodesic.
 Then, provided $\Delta$ has an $\cal R$-face,
there exists an $\cal R$-face $\Pi$ in $\Delta$ and disjoint $\varepsilon$-contiguity subdiagrams
$\Gamma_1,\dots,\Gamma_4$ (some of them may be absent) of $\Pi$ to $p_1$, $p_2$, $q_1$, $q_2$ respectively, such that
\begin{equation}  
(\Pi,\Gamma_1,p_1)+(\Pi,\Gamma_2,p_2)+(\Pi,\Gamma_3,q_1)+(\Pi,\Gamma_4,q_2) > 1-23\mu~.
\end{equation}
\end{lemma}

Collecting together the claims of Lemmas 6.7,7.4 and 7.5 from \cite{Olsh2} we obtain
\begin{lemma} \label{lemma6.7}
Suppose $G$ is a non-elementary hyperbolic group and $H'_1,\dots,H'_{k'}$ -- its non-elementary subgroups.
Then for any $\lambda >0$ there is $\mu_0 >0$ such that
for any $\mu \in (0,\mu_0]$ and $c \ge 0$ there is $\varepsilon \ge 0$ such that for any $N>0$ there exists
$\rho >0$ with the following property:

Let the symmetrized presentation (\ref{G_1}) satisfy the $C(\varepsilon,\mu,\lambda,c,\rho)$-condition.
Then the quotient $G_1$ (\ref{G_1}) is a non-elementary hyperbolic group and the images of the subgroups
$H'_1,\dots,H'_{k'}$ are non-elementary in $G_1$. Moreover, $W=1$ in $G_1$ if and only if $W=1$ in $G$
for every word $W$ with $\|W\| \le N$.
\end{lemma}

\section{Main Construction} \label{mainconstruction}
Assume, now, that we are in the conditions of Theorem \ref{mainthm}. The $\delta$-hyperbolic group $G$
is generated by the symmetrized set ${\cal A}=\{a'_1,\dots,a'_r\}$. Set $s=kr$ and
define $a_1,\dots,a_s$, $\H_1,\dots,\H_s$ as follows:
$$a_{ir+j}=a'_j,~\H_{ir+j}=H_{i+1} ~\mbox{ if }~ 1 \le j\le r,~0\le i \le k-1~, $$
i.e. $a_1=a'_1,\dots,a_r=a'_r,a_{r+1}=a'_1,a_{r+2}=a'_2, \dots$,
$\H_1=H_1,\dots,$ $\H_r=H_1,$ $\H_{r+1}=H_2,$ $\H_{r+2}=H_2,\dots$.

Since every $\H_i$ is a $G$-subgroup, we can find $b_i \in \H_i$ such that $a_ib_i^{-1} \in C_G\bigl(E(\H_i)\bigr)$
(such a choice is possible because $E(\H_i)=E(G)$ and $|\H_i:K(\H_i)|=|G:K(G)|$).

For each $i =1,2,\dots,s$, the subgroup $F_i=C_{\H_i}\bigl(E(\H_i)\bigr)$ has finite index in $\H_i$,
hence, $\Lambda(\H_i)=\Lambda(F_i)$ by parts $(c)$ and $(d)$ of lemma \ref{limitsets}. The set $Q^{-1}Q$ is quasiconvex
by lemma \ref{q-1q}, thus, according to the assumptions of theorem \ref{mainthm},
we can apply lemma \ref{findalpha} to find the points on the boundary $\partial G$:
$$\alpha_i \in \Lambda(F_i)\backslash \bigl(G\circ \Lambda(Q^{-1}Q)\bigr)$$
 and a sequence $\displaystyle \left(y^{(i)}_j\right)_{j \in \N} \subset F_i~ \mbox{ with } ~
\lim_{j \to \infty}y^{(i)}_j=\alpha_i,~ i=1,2,\dots,s.$

The set $\{y^{(i)}_j~|~j\in \N\}$ is infinite, therefore, the set $\{(y^{(i)}_j)^{-1}~|~j\in \N\}$
is also infinite, hence, by lemma \ref{limitsets} it has a limit point $\beta_i \in \partial G$. So, after passing to
a subsequence, we can assume that $$ \lim_{j \to \infty}(y^{(i)}_j)^{-1}=\beta_i,~i=1,2,\dots,s.$$

Using lemma \ref{lemma3.8} one can find an $\H_i$-suitable element $g_i\in {\H_i}^0$ for every $i=1,\dots,s$,
so that the elements $g_1,\dots,g_s$ are pairwise non-commensurable and on the Gromov boundary of the group
$G$ we have
\begin{equation} \label{g_i} \{g_i^{ \infty}, g_i^{-\infty}\} \cap \{\alpha_i,\beta_i\}=\emptyset,~i=1,2,\dots,s
\end{equation}
(recall that if for two elements of infinite order $g,h \in G$ we have $g^{\infty}=h^{\pm \infty}$ then
$g^m=h^l$ for some $m,l \in \Z \backslash \{0\}$ -- by lemma \ref{quasintersection}).

$\H_i$ is non-elementary, therefore $F_i \le \H_i$ is also non-elementary, $i=1,\dots,s$. Now we use lemma
\ref{findK} to obtain a non-elementary subgroup $K_i \le F_i$ such that
\begin{equation} \label{K}
\{ \alpha_i,~ g_i^{-\infty}, ~(b_ia_i^{-1})\circ g_i^{-\infty} \} \cap \Lambda(K_i)=\emptyset ~\mbox{ in }
  ~\partial G~,~i=1,2,\dots,s~.
\end{equation}


According to (\ref{g_i}) and (\ref{K}) we can use the claim of lemma \ref{boundedproduct} to show that
$$C_{01i}=\sup\left\{\left(\bigl(y^{(i)}_j\bigr)^{-1}|g_i^n\right)_{1_G}~:~j,n \in \N \right\}<\infty~, $$
$$C_{02i}=\sup\left\{\bigl(y^{(i)}_j|g_i^{-n}\bigr)_{1_G}~:~j,n \in \N \right\}<\infty~,$$
$$C_{03i}=\sup\left\{\bigl(x|y^{(i)}_j\bigr)_{1_G}~:~j \in \N, x \in K_i\right\}<\infty~,$$
$$C_{04i}=\sup\left\{\bigl(g_i^{-n}|x\bigr)_{1_G}~:~n \in \N, x \in K_i\right\}<\infty~\mbox{ and}$$
$$C_{05i}=\sup\left\{\bigl(b_ia_i^{-1}g_i^{-n}|x\bigr)_{1_G}~:~n \in \N,x \in K_i\right\}<\infty~$$
for each $i=1,2\dots,s$. Finally, define
$$C_0 = \max_{1\le i \le s}\left\{C_{01i},C_{02i},C_{03i}+|b_ia_i^{-1}|_G,C_{04i},C_{05i}+|b_ia_i^{-1}|_G\right\}
+14\delta~. $$

Denote $\bar \lambda = 1$, $\bar c = 0$,
\begin{equation} \label{C_1} C_1=12(C_0+\delta)+\bar c+1 \end{equation}
 Let
 \begin{equation} \label{lambda-c} \lambda=\bar \lambda/4 =1/4,~c=c(\bar \lambda,\bar c,C_0) \ge 0 ~\mbox{ be
the constant from lemma \ref{quasigeod}} \end{equation}
\begin{equation} \label{nu} \mbox{and let $\nu=\nu(\lambda,c)$ be the constant from the claim of lemma \ref{close}.}
\end{equation}

By the conditions of the theorem \ref{mainthm}, the subset $Q$ is $\eta$-quasiconvex for some $\eta \ge 0$,
\begin{equation} \label{varkappa} \mbox{let $\varkappa$ be the length of a shortest element from $Q$.}
\end{equation}

Set
\begin{equation} \label{A} A=\{g \in G~|~|g|_G\le 3\delta+\nu+\eta+\varkappa \}
\end{equation} then $card(A)<\infty$. By the construction of $\alpha_i$,
$$\alpha_i \notin \bigcup_{g \in A} g\circ\Lambda(Q^{-1}Q)=\Lambda\left(
\bigcup_{g \in A} gQ^{-1}Q\right)=\Lambda(AQ^{-1}Q)~.$$

Hence, according to the lemma \ref{boundedproduct}, one can define
\begin{equation} \label{C_3} C_3 = \max_{1\le i \le s} \sup \{ \bigl(y^{(i)}_j|x  \bigr)_{1_G}~|
~ j\in \N, x \in AQ^{-1}Q\}~. \end{equation}

For every $i\in \{1,2,\dots,s\}$, ~$|y^{(i)}_j|_G \to \infty$ as $j\to \infty$, and the intersection
$\{y^{(i)}_j~|~j\in \N\} \cap E(g_i)$
is finite by (\ref{g_i}), therefore for some $j_0$ (depending on $i$), after setting $y_i=y^{(i)}_{j_0} \in F_i$, we will have
\begin{equation} \label{y_i} |y_i|_G > 3\delta+\nu+\eta+\varkappa + C_1+2C_3 ~\mbox{ and }~y_i \in
C_{\H_i}\bigl(E(\H_i)\bigr) \backslash E(g_i)~.
\end{equation}

Applying lemmas \ref{y-mod} and \ref{y-modnon-comm} we can find $n \in \N$ such that
the elements
$$w_i=y_ig_i^n$$ have infinite order, are $\H_i$-suitable and pairwise non-commensurable when \\ $1\le i \le s$.
Obviously, in addition, we can demand that ~$|g_i^n|_G,|w_i|_G>C_1$ for every $i$.

The subgroups $K_i \le G$ are non-elementary, hence we can find elements \\
$c_i \in K_i \le C_{\H_i}\bigl(E(\H_i)\bigr)$ for which
\begin{equation} \label{x_i0} x_{i0}\stackrel{def}{=}a_ib_i^{-1}c_i \in C_{G}\bigl(E(\H_i)\bigr) \backslash E(w_i)~
\mbox{ and }~ |x_{i0}|_G > C_1 ~\mbox{ for }~ i=1,\dots,s~.\end{equation}
Note that $b_i^{-1}c_i \in \H_i$ for each $i$.

Let $w_i \in \hat H_i^0$, $x_{i0} \in G$, $x_{i1},\dots,x_{il} \in \hat H_i$ be represented by
words $W_i$, $X_{i0}$, $X_{i1},\dots,X_{il}$ over the alphabet $\cal A$ of minimal length, $i=1,\dots,s$. As the system
${\cal R}={\cal R}_{s,l,m}(W_1,\dots,W_s,X_{10},\dots,X_{sl},m)$ consider all cyclic permutations of $R_i^{\pm 1}$
where
\begin{equation} \label{R_i} R_i \equiv X_{i0}W_i^mX_{i1}W_i^m \dots X_{il}W_i^m~,~~ i=1,2,\dots,s~.\end{equation}

\begin{lemma} \label{R-quasiged} Fix $i \in \{1,\dots,s\}$. Let $\lambda,c$ be the constants defined in (\ref{lambda-c})
 and the elements $x_{ij}$, $j=1,\dots,l$,  satisfy the properties:
$$x_{ij} \in K_i~\mbox {and } |x_{ij}|_G>C_1~\mbox{ for every $j=1,\dots,l$}$$
($K_i\le \H_i$ and $C_1 >0$ are as above).
Then for any $n \in \N$, any path $q$ in \ga labelled by a word $R_i^{\pm n}$ is
$(\lambda,c)$-quasigeodesic.
\end{lemma}

\underline{Proof.} It is enough to prove this lemma for the case when $q$ is labelled by $R_i^n$ because if
$lab(q) \equiv R_i^{-n}$ then $lab(q^{-1}) \equiv R_i^{n}$ and if $q^{-1}$ is $(\lambda,c)$-quasigeodesic then
so is $q$.

For convenience, assume that $i=1$. Since left translations are isometries,
we can suppose that $q_-=1_G$. Set $t=(m+1)(l+1)$. The path $q$ is a broken line
$[z_0,z_1,\dots,z_{nt}]$ in the Cayley graph \ga where

$z_0=q_-=1_G$, $z_1=x_{10}$, $z_2=x_{10}w_1$, $z_3=x_{10}w_1^2$,

$\dots\dots$

$z_{m+1}=x_{10}w_1^m$, $z_{m+2}=x_{10}w_1^mx_{11}$, $z_{m+3}=x_{10}w_1^mx_{11}w_1$,

$\dots\dots$

$z_{t}=x_{10}w_1^mx_{11}w_1^m \dots x_{1l}w_1^m$, $z_{t+1}=x_{10}w_1^mx_{11}w_1^m \dots x_{1l}w_1^mx_{10}$,

$\dots\dots$

$z_{nt}=x_{10}w_1^m \dots x_{1l}w_1^mx_{10}w_1^m \dots x_{1l}w_1^m \dots  x_{10}w_1^m \dots x_{1l}w_1^m=q_+$.

By construction, $\|[z_{j-1},z_j]\|>C_1$, $j=1,2,\dots,nt$.
In order to apply lemma \ref{quasigeod} to the path $q$ it remains to verify that $(z_{j-1}|z_{j+1})_{z_j} \le C_0$ for
every $l=1,\dots,nt-1$. There are several types of Gromov products that appear when $j$ changes from $1$ to $nt-1$.
Below we compute them in the order of their occurrence.

{\bf Type I.} $\displaystyle (z_0|z_2)_{z_1}=(1_G|x_{10}w_1)_{x_{10}}=(x_{10}^{-1}|w_1)_{1_G}.$
Recall that $w_1=y_1g_1^n$. By the Gromov's definition of a hyperbolic space,
$$ (x_{10}^{-1}|y_1)_{1_G} \ge \min\{(x_{10}^{-1}|w_1)_{1_G},(y_1|w_1)_{1_G}\} -\delta~.$$

Now, we observe that $$(x_{10}^{-1}|y_1)_{1_G} = (c_1^{-1}b_1a_1^{-1}|y_1)_{1_G} \le
(c_1^{-1}|y_1)_{1_G}+|b_1a_1^{-1}|_G \le C_{031}+|b_1a_1^{-1}|_G.$$
Hence, \begin{equation} \label{gromprod1} \min\{(x_{10}^{-1}|w_1)_{1_G},(y_1|w_1)_{1_G}\} \le C_{031}+|b_1a_1^{-1}|_G
+\delta \le C_0~. \end{equation}

From the geodesic triangle $1_Gy_1w_1$ in \ga we obtain
$$(y_1|w_1)_{1_G}=|y_1|_G-(1_G|w_1)_{y_1} = |y_1|_G-(y_1^{-1}|g_1^n)_{1_G} \ge C_1-C_0>C_0~.$$

Combining the latter inequality with (\ref{gromprod1}) we achieve
$$(z_0|z_2)_{z_1}=(x_{10}^{-1}|w_1)_{1_G} \le C_0.$$

{\bf Type II.} $\displaystyle (z_1|z_3)_{z_2}=(x_{10}|x_{10}w_1^2)_{x_{10}w_1}=(w_1^{-1}|w_1)_{1_G}.$

Again, applying the definition of hyperbolicity twice, we obtain
$$(g_1^{-n}|y_1)_{1_G} \ge \min\{(g_1^{-n}|w_1)_{1_G},(y_1|w_1)_{1_G}\}-\delta \ge $$
$$ \ge \min\{(w_1^{-1}|w_1)_{1_G},(g_1^{-n}|w_1^{-1})_{1_G},(y_1|w_1)_{1_G}\}-2\delta.$$

By construction, $(g_1^{-n}|y_1)_{1_G} \le C_{021} \le c_0-2\delta$. As we showed above, $(y_1|w_1)_{1_G} > C_0$.

Considering the geodesic triangle $1_Gg_1^{-n}w_1^{-1}$ we get
$$(g_1^{-n}|w_1^{-1})_{1_G}=|g^{-n}|_G-(1_G|w_1^{-1})_{g^{-n}}= |g^{n}|_G-(g^n|y_1^{-1})_{1_G}\ge C_1-C_0>C_0.$$

So, combining these inequalities, we achieve
$$(z_1|z_3)_{z_2}=(w_1^{-1}|w_1)_{1_G} \le C_0~.$$

{\bf Type III.} $\displaystyle (z_m|z_{m+2})_{z_{m+1}}=(w_1^{-1}|x_{11})_{1_G}\le C_0.$

{\bf Type IV.} $\displaystyle (z_{m+1}|z_{m+3})_{z_{m+2}}=(x_{11}^{-1}|w_1)_{1_G}\le C_0.$ \\
These two inequalities are proved in the same way as we proved the inequality for Type I
(the proofs even easier since $x_{11} \in K_1$). \\
The last possibility is

{\bf Type V.} $\displaystyle (z_{t-1}|z_{t+1})_{z_{t}}=(w_1^{-1}|x_{10})_{1_G}.$

As before, we have
$$(g_1^{-n}|x_{10})_{1_G} \ge \min \{(w_1^{-1}|x_{10})_{1_G},(w_1^{-1}|g_1^{-n})_{1_G}\} -\delta~.$$
$$(g_1^{-n}|x_{10})_{1_G} =(g_1^{-n}|a_1b_1^{-1}c_1)_{1_G}=(b_1a_1^{-1}g_1^{-n}|c_1)_{b_1a_1^{-1}}\le $$
$$\le (b_1a_1^{-1}g_1^{-n}|c_1)_{1_G}+|b_1a_1^{-1}|_G \le C_{051}+|b_1a_1^{-1}|_G \le C_0 -\delta~.$$

We showed while considering the Type II, that $(g_1^{-n}|w_1^{-1})_{1_G} > C_0$.
Therefore,
$$(z_{t-1}|z_{t+1})_{z_{t}}=(w_1^{-1}|x_{10})_{1_G} \le C_0.$$

It is easy to see that for arbitrary $j \in \{1,2,\dots,nt-1\}$ the Gromov product $(z_{j-1}|z_{j+1})_{z_j}$
is equal to a Gromov product of one of the Types I-V, thus it is not larger than $C_0$.

Therefore, recalling that the constant $C_1$ was defined by formula (\ref{C_1}), we can use the lemma \ref{quasigeod}
to show that the path $q$ is $(\lambda,c)$-quasigeodesic, where $\lambda>0$ and $c\ge 0$ are defined in (\ref{lambda-c}).
Q.e.d. $\square$

\vspace{.15cm}
Below we have an analog of the lemma \ref{lemma4.2} needed for our proof:

\begin{lemma} \label{smallcancell-mod}
Suppose $W_1,\dots,W_s$, $X_{10},\dots,X_{s0}$ and $\lambda >0$, $c\ge 0$ are the words and the constants
defined above. Then
for any $\mu >0$ there are $l \in \N$ and words
$X_{11},\dots,X_{sl}$ ($X_{ij}$ represents an element $x_{ij} \in \H_i$, $j=1,\dots,l$, $i=1,\dots,s$)
such that for any $\varepsilon \ge 0$, $\rho>0$ there is $m_0 \in \N$
such that the system ${\cal R}_{s,l,m}$ (\ref{R_i})
satisfies $C(\varepsilon,\mu,\lambda,c,\rho)$ and $C_1(\varepsilon,\mu,\lambda,c,\rho)$-conditions if $m \ge m_0$.
\end{lemma}

\underline{Proof.} By lemma \ref{lemma4.2} there exist $\lambda' >0$ such that for any $\mu>0$ there are $l \in \N$
and $c' \ge 0$ such that for any $\varepsilon \ge 0$, $\rho>0$ there are $m_0 \in \N$ and words
$X_{11},\dots,X_{sl}$ such that the system ${\cal R}_{s,l,m}$ satisfies the generalized small cancellation
conditions $C(\varepsilon,\mu,\lambda',c',\rho)$ and
$C_1(\varepsilon,\mu,\lambda',c',\rho)$ if $m \ge m_0$.

According to the remark \ref{remafter4.2} after the formulation of lemma \ref{lemma4.2},
lemma \ref{lemma3.7} and lemma \ref{lemma3.8}, the elements $x_{i1},\dots,x_{is}$ can be
chosen right after $l$, inside of the subgroup $K_i$, with an additional property $|x_{ij}|_G>C_1$
(the constant $C_1$ was defined in (\ref{C_1})) for every $j=1,\dots,l$, $i=1,\dots,s$.

Consider any word $R \in {\cal R}_{s,l,m}$. By definition, $R$ is a subword of a word $R_i^{\pm2}$
for some $i \in \{1,\dots,s\}$. By lemma
\ref{R-quasiged} the word $R_i^{\pm2}$ is $(\lambda,c)$-quasigeodesic (where $\lambda$
and $c$ are defined in (\ref{lambda-c})), hence so is $R$.
Taking into account remark \ref{indepoflambdac}, we achieve
that the system ${\cal R}_{s,l,m}$ satisfies the conditions $C(\varepsilon,\mu,\lambda,c,\rho)$ and
$C_1(\varepsilon,\mu,\lambda,c,\rho)$ if $m \ge m_0$. $\square$

\begin{lemma} \label{shortcontig} Let ${\cal R}={\cal R}_{s,l,m}(W_1,\dots,W_s,X_{10},\dots,X_{sl},m)$ be the system
of additional
relations and $\lambda >0$, $c \ge 0$ be the constants defined above. Then for any $\varepsilon \ge 0$ and
$\xi >0$ there exists $m_1 \in \N$ such that for any $m \ge m_1$ the following property holds:

Suppose $\Delta$ is a diagram over the presentation (\ref{G_1}) and $q$ -- a subpath of $\partial \Delta$
such that the corresponding path $q'$ in the Cayley graph \ga of the group $G$ with the same label as $q$
is geodesic (in other words, $\|q\|=|elem(q)|_G$) and $elem(q) \in Q$ in $G$. Then for arbitrary
$\cal R$-face $\Pi$ of $\Delta$ and an $\varepsilon$-contiguity subdiagram $\Gamma$ between $\Pi$ and $q$ we have
$$(\Pi,\Gamma,q) \le \xi~.$$
\end{lemma}

\underline{Proof.} Let $\partial \Gamma = p_1q_1p_2q_2$ where $q_1$, $q_2$ are subpaths of $\partial \Pi$ and $q$
correspondingly and $\|p_1\|,\|p_2\| \le \varepsilon$. Take arbitrary $\xi >0$.
Obviously, from the definition (\ref{R_i}),
 there is $m_1 \in \N$ such that for any $m \ge m_1$ the inequality $\|q_1\|/\|\partial \Pi\| > \xi$
implies that the path $q_1$ has a subpath $o$ labelled by the word $W_i^{\pm 1}$ for some $i\in \{1,\dots,s\}$ and,
 moreover, the subpaths $o_1$, $o_2$ of $q_1$ (with $(o_1)_-=(q_1)_-$, $(o_1)_+=o_-$, $(o_2)_-=o_+$,
$(o_2)_-=(q_1)_+$) satisfy
\begin{equation} \label{o_j} \|o_j\| > (\varepsilon +c+\nu)/\lambda,~ j=1,2~,\end{equation}
($\nu$ is chosen according to (\ref{nu})).

We are going to obtain a contradiction with the definitions of elements $w_i$ and $y_i$.

Since the diagram $\Gamma$ contains only $0$-faces (i.e. it is a diagram over the group $G$), we can consider
the corresponding picture in \ga with
a geodesic path $q'$ starting at $1_G$ (its subpath $q_2'$),
$(\lambda,c)$-quasigeodesic path $q_1'$ (its subpaths $o'$, $o_1'$, $o_2'$) and paths $p_1'$, $p_2'$
of lengths at most $\varepsilon$ with $(p_1')_-=(q_2')_+$, $(p_1')_+=(q_1')_-$, $(p_2')_-=(q_1')_+$,
$(p_2')_+=(q_2')_-$ (i.e. for every path $r$ from $\Delta$ we construct a corresponding path $r'$
in \ga with the same label; see Figure 1).

Pick any $z \in Q$  with $|z|_G=\varkappa$ (the constant $\varkappa$ was defined in (\ref{varkappa})).
Then $q'_+=elem(q')=elem(q) \in Q$. Hence, since the triangles are $\delta$-slim, one obtains
$$q' \subset {\cal O}_{\delta}([1,z] \cup [z,elem(q')]) \subset {\cal O}_{\delta +\varkappa}([z,elem(q')])
\subset {\cal O}_{\delta +\varkappa+\eta}(Q)~.$$

Denote $u=(q_1')_-$, $v=(q_1')_+$. Then $\displaystyle u,v \subset {\cal O}_{\varepsilon}(q_2')$.

\begin{figure}[!ht]
   \begin{center}
    \input{pic1.tex}

   Figure 1
   \end{center}
\end{figure}

Now, since $$ o' \subset q_1' \subset {\cal O}_{\nu} ([u,v])~,$$ using (\ref{o_j}) and lemma \ref{quadrangle}
we obtain $$ o' \subset {\cal O}_{\nu +2\delta}(q') \subset {\cal O}_{3\delta +\nu+\varkappa+\eta}(Q)~.$$

Recall that $lab(o')=W_i^{\pm 1}$ by construction. So, if $lab(o')=W_i$, define the points $f=o'_-$,
$g=o'_+$ and if $lab(o')=W_i^{-1}$, define $g=o'_-$ and $f=o'_+$. Thus there are elements $h_1,h_2 \in Q$
such that
$$d(f,h_1) = |f^{-1}h_1|_G \le 3\delta +\nu+\varkappa+\eta, ~d(g,h_2) \le 3\delta +\nu+\varkappa+\eta.$$

By the definition of $w_i$, we have $fy_ig_i^n=g$. According to Gromov's definition of a
hyperbolic metric space, we achieve
$$(h_2|fy_i)_f \ge \min \{(h_2|g)_f,(g|fy_i)_f\} -\delta~.$$
Observe that $(h_2|fy_i)_f=(f^{-1}h_2|y_i)_{1_G}$ and $$x=f^{-1}h_2=(f^{-1}h_1)h_1^{-1}h_2\in AQ^{-1}Q$$
(the set $A$ was defined in (\ref{A})).

$$(g|fy_i)_f=\|[f,fy_i]\|-(f|g)_{fy_i} = |y_i|_G - (y_i^{-1}|g_i^n)_{1_G} \ge$$ $$ \ge |y_i|_G-C_0>
5\delta +\nu+\varkappa+\eta+2C_3~.$$
(Here we used that $C_1-C_0 > 2\delta$.) Note that $d(f,g)\ge (g|fy_i)_f$, hence
$$(h_2|g)_f \ge \frac12\bigl(d(f,g)-d(g,h_2)\bigr)\ge \frac12\bigl((g|fy_i)_f-(3\delta +\nu+\varkappa+\eta)\bigr)>
C_3+\delta~.$$

Combining the above formulas, we finally obtain
$$(h_2|fy_i)_f = (x|y_i)_{1_G} > C_3$$
contradicting to the definition (\ref{C_3}) of $C_3$. Therefore, $\|q_1\|/\|\partial \Pi\| = (\Pi,\Gamma,q) \le \xi$.

The lemma is proved. $\square$

\vspace{.15cm}

\section{Proof of Theorem \ref{mainthm}}
\underline{Proof.} The group $G_1$ is generated by
$\phi({\cal A})$, so let $|x|_{G_1}$ be the corresponding length function for elements $x \in G_1$, and let
$d_1(\cdot,\cdot)$ be the corresponding metric on the Cayley graph of the group $G_1$.
Sometimes it will be convenient for us to identify $\cal A$ and $\phi({\cal A})$ for $G_1$, so $\Gamma(G_1,{\cal A})$
will be the Cayley graph of $G_1$.
Since $\phi$ is a
homomorphism, from  the definition of the word metric it follows that
\begin{equation} \label{metrics} \forall~x,y\in G~~d_1\bigl(\phi(x),\phi(y)\bigr) \le d(x,y)~.
\end{equation}


Define the elements $a_1,\dots,a_s \in G$ and the subgroups $\H_1,\dots,\H_s$ as we did in the beginning of
section \ref{mainconstruction}. After that construct the elements $g_i$, $y_i$, $w_i$ and $x_{i0}$, $i=1,2,\dots,s$,
as described in that section. Then we can find the constants $\lambda>0$ and $c\ge 0$ according to (\ref{lambda-c}).

Suppose that $W_i$, $X_{i0},\dots,X_{il}$ are shortest words in the alphabet $\cal A$ representing
$w_i$, $x_{i0},\dots,x_{il}$, $i=1,\dots,s$.
As the system of additional relations, consider the set
$${\cal R}={\cal R}_{s,l,m}(W_1,\dots,W_s,X_{10},\dots,X_{sl},m)$$
of all cyclic permutations of $R_i^{\pm 1}$, $i=1,\dots,s$, established in (\ref{R_i}).

Define the group $G_1$ according to (\ref{G_1}), thus, $G_1 \cong G/\langle {\cal R}^G \rangle$. Let $\phi$
be the natural epimorphism from $G$ to $G_1$.

By lemma \ref{smallcancell-mod} one can find $l,m_0 \in \N$ and elements $x_{ij}\in \H_i$, $j=1,\dots,l$,
$i=1,\dots,s$, such that the group $G_1$ satisfies
all of the conditions of lemmas \ref{lemma6.6} and \ref{lemma6.7} if $m \ge m_0$. Therefore we obtain
the parts 1) and 8) of the theorem \ref{mainthm}.

It is easy to see that the relation $R_i$ implies $\phi(a_iz_i)=1$ in $G_1$ for some $z_i \in \H_i$,
hence $\phi(a_i) \in \phi(\H_i)$ for $i=1,\dots,s$.

Due to the choice of $a_1,\dots,a_s$ and $\H_1,\dots,\H_s$ we obtain $\phi({\cal A}) \subset \phi(H_j)$
for every $j \in \{1,2,\dots,k\}$. Consequently, $G_1=\phi(H_j)$, $j=1,\dots,k$, so the part 3) of the theorem
is proved.

Let us now prove the property 2). Let $\mu_0 >0$, $\varepsilon \ge 0$ be chosen according to lemma \ref{lemma6.6}.
Since we can take any $\mu$ inside of the interval $(0,\mu_0]$ we can also demand that
~$1/(\lambda+1)<1-23\mu $. Choose $\xi >0$ in such a way that
\begin{equation} \label{xi} \frac1{\lambda+1}<1-23\mu-2\xi ~.
\end{equation}

Denote $\theta = 1-23\mu-2\xi>0$. Then (\ref{xi}) implies that $(\lambda+1)\theta-1>0$.
Set $L_0=\min \{\|R\|~|~R\in {\cal R}\}$. Evidently, $L_0$ depends on $m$ and
 there exists $m_2 \in \N$ such that for any $m \ge m_2$
\begin{equation} \label{m_2} \bigl((\lambda+1)\theta-1\bigr)L_0>c+4\varepsilon~.
\end{equation}

Now, let's apply the statement lemma \ref{shortcontig} to find $m_1=m_1(\varepsilon,\xi) \in \N$.

By taking any $m \ge \max\{m_0,m_1,m_2\}$ we can further assume that the claims of lemmas \ref{lemma6.6}
and \ref{shortcontig} hold together with the inequality (\ref{m_2}).

Consider arbitrary elements $u,v \in Q$. We need to show that
$d(u,v)=d_1\bigl(\phi(u),\phi(v)\bigr)$.

Observe that, by definition, $d(u,v)=|u^{-1}v|_G$, $d_1(\phi(u),\phi(v))=
|\phi(u^{-1}v)|_{G_1}$. Obviously, $|u^{-1}v|_G \ge |\phi(u^{-1}v)|_{G_1}$, so assume, by contradiction,
that \begin{equation} \label{assumption} |u^{-1}v|_G > |\phi(u^{-1}v)|_{G_1}~. \end{equation}

Thus, if $U,V$ are shortest words representing $u,v$ in $G$, there is a word $Z$ such that
$U^{-1}V=Z$ in the $G_1$ but not in $G$ ($Z$ is a word of minimal length representing
the element $\phi(u^{-1}v)$ in $G_1$).

Consider a reduced circular diagram $\Delta$ over $G_1$ whose boundary is labelled by the word
$U^{-1}VZ^{-1}$. Let $q^1$, $q^2$, $q^3$ be the (geodesic) sections of the boundary $\partial \Delta$
labelled by the words $U$, $V$, $Z$ respectively.

This diagram must contain at least one $\cal R$-face since $U^{-1}VZ^{-1}\neq 1$
in $G$. Therefore, by lemma \ref{lemma6.6} there exists an $\cal R$-face $\Pi$ in $\Delta$ and
$\varepsilon$-contiguity subdiagrams $\Gamma_1$,$\Gamma_2$,$\Gamma_3$ between $\Pi$ and the sections $q^1$,$q^2$,$q^3$
(for our convenience, for each of the sections $q^j$ we can choose a corresponding orientation of $\partial \Pi$,
$j=1,2,3$) satisfying $$(\Pi,\Gamma_1,q^1)+(\Pi,\Gamma_2,q^2)+(\Pi,\Gamma_3,q^3) > 1-23\mu~.$$

Since $elem(q^1)=u \in Q$, $elem(q^2)=v \in Q$ and $m\ge m_1$, we have
$(\Pi,\Gamma_1,q^1) \le \xi$ and $(\Pi,\Gamma_2,q^2) \le \xi$. Hence,
\begin{equation} \label{q^3} (\Pi,\Gamma_3,q^3) > 1-23\mu-2\xi=\theta~. \end{equation}

Now we are going to obtain a contradiction with the choice of $Z$.
\\ Let $\partial(\Gamma_3)=p_1r_1p_2o_2$ where $\partial \Pi=r_1r_2$, $q^3=o_1o_2o_3$, $\|p_1\|,\|p_2\|\le \varepsilon$
(Figure 2).

  \begin{figure}[!ht]
   \begin{center}
    \input{pic2.tex}

   Figure 2
   \end{center}
  \end{figure}

Let $L$ denote the length of $\partial \Pi$. (\ref{q^3}) implies
\begin{equation}  \label{r_1} \|r_1\| > \theta L~,~~ \|r_2\|=L-\|r_1\|<(1-\theta)L~.
\end{equation}

Now, since $\Gamma_3$ is a diagram over the group $G$, the equality
$$elem(o_2^{-1})=elem(p_1)elem(r_1)elem(p_2)$$
holds in $G$. The path $q^3$ is geodesic, therefore, its subpath $o_2$ is also geodesic, thus,
 $$\|o_2\|=\|o_2^{-1}\|=|elem(o_2^{-1})|_G \ge |elem(r_1)|_G - |elem(p_1)|_G-|elem(p_2)|_G~,~\mbox{hence}$$
$$ \|o_2\| \ge  |elem(r_1)|_G- 2\varepsilon~.$$

The path $r_1$ is $(\lambda,c)$-quasigeodesic as a subpath of the face contour $\partial \Pi$,
therefore $|elem(r_1)|_G\ge \lambda\|r_1\|-c$. Combining the last two inequalities with (\ref{r_1}) we obtain
 \begin{equation} \label{o_2} \|o_2\| \ge \lambda \theta L -c- 2\varepsilon~.\end{equation}

Consider the subdiagram $\Omega$ of $\Delta$ bounded by the closed path $p_2^{-1}r_2p_1^{-1}o_2^{-1}$.
It gives us the following equality in the group $G_1$:
$$\label{elemo_2} elem(o_2)=elem(p_2^{-1}) \cdot elem(r_2) \cdot elem(p_1^{-1})~.$$ Thus,

$$ \|o_2\| = |elem(o_2)|_{G_1} \le |elem(p_2^{-1})|_{G_1} + |elem(r_2)|_{G_1} + |elem(p_1^{-1})|_{G_1} \le $$
$$\le \|r_2\|+2\varepsilon \le (1-\theta)L+2\varepsilon ~. $$
Comparing the latter inequality with (\ref{o_2}) we get
$$\lambda \theta L -c- 2\varepsilon \le (1-\theta)L+2\varepsilon~.$$
Or, equivalently, $$\bigl((\lambda+1)\theta-1\bigr)L\le c+4\varepsilon~.$$
Since $L\ge L_0$ this contradicts to the inequality (\ref{m_2}).

Therefore, the assumption (\ref{assumption}) was incorrect and $d(u,v)=d_1\bigl(\phi(u),\phi(v)\bigr)$
for arbitrary $u,v \in Q$. Thus $\phi|_Q$ is an isometry.

By 1) $G_1$ is $\delta_1$-hyperbolic for some $\delta_1\ge 0$. Take any $\omega$-quasiconvex (in $G$) subset
$S \subseteq Q$. Let's show that $\phi(S) \subset G_1$ is $(\omega+\delta_1)$-quasiconvex.

Consider arbitrary two elements $u,v \in S$ and let $p$ be a geodesic path in \ga connecting them. Then
$$p\subset {\cal O}_{\omega}(S)~\mbox{ in } \mbox{\ga}~.$$
Let $p_1$ be the path in $\Gamma(G_1,{\cal A})$ starting at $\phi(u)$ with the same label as $p$.
Then $(p_1)_+=\phi(v)$ (this is equivalent to the equality $\phi(u)\cdot elem(p_1)=\phi(v)$ which follows from
$u  \cdot elem(p)=v$). Now, since $\phi$ is an isometry between $S$ and $\phi(S)$,
$$\|p_1\|=\|p\|=d(u,v)=d_1\bigl(\phi(u),\phi(v)\bigr)~.$$
Therefore, $p_1$ is a geodesic path between $\phi(u)$ and $\phi(v)$ in $\Gamma(G_1,{\cal A})$. (\ref{metrics})
implies $$p_1 \subset {\cal O}_{\omega}\bigl(\phi(S)\bigr)~\mbox{ in } \Gamma(G_1,{\cal A})~.$$
The space $\Gamma(G_1,{\cal A})$ is $\delta_1$-hyperbolic, hence for any geodesic path $q$ between
$\phi(u)$ and $\phi(v)$ we have $\displaystyle q \subset {\cal O}_{\delta_1}(p_1)$. Consequently,
$$q \subset {\cal O}_{\omega+\delta_1}\bigl(\phi(S)\bigr)~\mbox{ in } \Gamma(G_1,{\cal A})~.$$
The proof of the part 2) is complete.

For the case when $Q$ is a finite subset, the proofs in \cite[Thm. 3]{Olsh2} of the properties corresponding to
$4),5)$ from our theorem \ref{mainthm} were based on the lemma \ref{lemma8.1}, general properties of hyperbolic groups
and the fact that in a diagram over (\ref{G_1}) with labels of boundary contours representing elements of $Q$ in $G$
there can not exist any "long" $\varepsilon$-contiguity of an $\cal R$-face to a boundary contour.
The same fact is true in our case by lemma \ref{shortcontig} (after an appropriate choice the parameters like
in the proof of the property 2)). So, for the proofs $4),5)$ the reader is referred to \cite[Thm. 3]{Olsh2}.

Properties $6)$ and $7)$ do not depend on $Q$, thus they can be proved in the same way as they were proved
in \cite[Thm. 2]{Olsh2} (we can always add a finite subset to $Q$: it will stay quasiconvex
and the formula $(*)$ will continue to hold).

Finally, let's derive the property 9). By lemma  \ref{lemma3.8}, we can choose a $G$-suitable element $g \in G$.
Then, by definition, $T(g)=E(G)$. Denote $S=\langle g \rangle_\infty$ -- a quasiconvex subgroup of the group $G$.
Then for any $h\in G$ $$|H_i:(H_i\cap hSh^{-1})|=\infty$$ since $H_i$ is non-elementary for every $i=1,\dots,k$.
Hence, according to theorem \ref{qcex}, $H_i$ is not contained inside of $P_1S^{-1}SP_2$ for arbitrary finite subsets
$P_1,P_2 \subset G$. Now we can apply lemmas \ref{q-cunion,product} and \ref{smallunion} to the union
$$Q'=Q \cup S=Q \cup \langle g \rangle_\infty$$ to show
that all the requirements of theorem \ref{mainthm} will remain satisfied if one substitutes $Q$ by $Q'$ in it.
Since the properties 1)-8) were already proved, we can further use them for the elements of $Q'$.
Therefore, $ker(\phi) \cap Q' =\{1_G\}$, implying that $\phi(g)$ has infinite order in $G_1$.

Consider arbitrary $x \in E(G_1)$. Then, in particular, $x \in E\bigl(\phi(g)\bigr)$.
By definition, there exists $n \in \N$ such that $x(\phi(g))^n x^{-1}= (\phi(g))^{\pm n}$.

If $x\phi(g)^n x^{-1}= \phi(g)^{-n}$ then by the part 4) the elements $g^n, g^{-n} \in Q'$ must be conjugate in $G$
which fails because $E(g)=E^{+}(g)$. Hence, $x(\phi(g))^n x^{-1}= (\phi(g))^{n}$, i.e.
$x \in C_{G_1}\bigl(\phi(g^n)\bigr)$.

Since $g^n \in Q'$, one can apply the part 5) to find $y \in C_G(g^n)$ with \\ $\phi(y)=x$.
$g \in G$ is $G$-suitable, therefore $C_G(g^n) \le E(g) = T(g) \times \langle g \rangle$. $G_1$ is non-elementary,
therefore the subgroup $E(G_1)\le G_1$ is finite, thus $x$ has a finite order in $G_1$. It follows that $y$ has a
finite order in $G$, because, otherwise, we would get $y^{l_1}=g^{l_2}$ for some $l_1,l_2 \in \Z \backslash \{0\}$ and
$x^{l_1}=\phi(y^{l_1})=\phi(g^{l_2})$ where $\phi(g^{l_2})$ has an infinite order in $G_1$. Consequently,
$y \in T(g) = E(G)$ and $$x=\phi(y) \in \phi\bigl(E(G)\bigr)~.$$

The proof of the theorem is finished. $\square$

\section{Constructing Simple Quotients}

\begin{lemma} \label {N-K} Suppose $N$ is an infinite normal subgroup of a hyperbolic group $G$ and $K$ is a quasiconvex
subgroup of $G$ such that $|G:K|=\infty$. Then for arbitrary $h \in G$, $|N:(N\cap hKh^{-1})|=\infty$.
\end{lemma}

\underline{Proof.} Since a conjugate to a quasiconvex subgroup of infinite index is again a quasiconvex subgroup of
infinite index, it is enough to consider the case when $h=1_G$.
Assume, by the contrary, that $|N:(N\cap K)|<\infty$. Then there exist elements $h_1,\dots,h_n \in N$
such that $N \subseteq Kh_1 \cup \dots \cup Kh_n$. Applying lemmas \ref{limitofnormal} and \ref{limitsets}  we achieve
$$\Lambda(G) = \Lambda(N) \subseteq \Lambda(Kh_1 \cup \dots \cup Kh_n)=\Lambda(K)~.$$

Hence, by lemma \ref{tame-qc}, $G \subset K\cdot P=\bigcup_{p \in P} Kp$~ for some finite subset $P$ of $G$, which
implies that $|G:K|< \infty$ -- a contradiction to our conditions. $\square$


\begin{lemma} \label{kernelofaction} If $H$ is a non-elementary subgroup of a hyperbolic group $G$ then $E(H)$
coincides with the subgroup $\displaystyle K =\bigcap_{\alpha \in \Lambda(H)}St_G(\{\alpha\}) \le G$.
\end{lemma}

\underline{Proof.} Indeed, $$K \subseteq
\bigcap_{h \in H^0} St_G(\{h^\infty\})=\bigcap_{h \in H^0} E^+(h) \subseteq E(H)~,$$
i.e. $K \subseteq E(H)$. Now, since $E(H)$ is a finite subgroup normalized by $H$, for every $x \in E(H)$
and $\alpha \in \Lambda(H)$ we can find a sequence of elements $y_i \in H$, $i \in \N$, with
$\lim_{i \to \infty} y_i = \alpha$. Denote $Y=\{y_i~|~i \in \N\} \subset H$. Then $\Lambda(Y)=\{\alpha\}$,
$xY \subset Y\cdot E(H)$ and by lemma \ref{limitsets} we achieve
$$x \circ \{\alpha\} =x \circ \Lambda(Y) = \Lambda(xY) \subset  \Lambda\bigl(Y E(H)\bigr)=$$
$$=\Lambda(\bigcup_{z \in E(H)} Yz) = \bigcup_{z \in E(H)}\Lambda(Yz) = \bigcup_{z \in E(H)}\Lambda(Y)=\{\alpha\}.$$
Thus, $x\in K$ which implies $E(H) \subseteq K$. $\square$

\begin{lemma} \label{elementalizereq} Assume that $A$ is a non-elementary normal subgroup of a subgroup
$H$ in a hyperbolic group $G$. Then $E(A)=E(H)$.
\end{lemma}

\underline{Proof.} 
According to lemma \ref{kernelofaction} and lemma \ref{limitofnormal} we have
$$E(A)= \bigcap_{\alpha \in \Lambda(A)} St_G(\{\alpha\})=\bigcap_{\alpha \in \Lambda(H)} St_G(\{\alpha\})=E(H)~.$$
$\square$

\underline{Proof of corollary \ref{simpleinductive}.} First, since $E(G)$ is the maximal finite normal subgroup of $G$, we can
consider the quotient $\hat G=G/E(G)$. Obviously, the natural homomorphism $\psi: G \to \hat G$ is a quasiisometry between
$G$ and $\hat G$, therefore $\hat G$ is a non-elementary hyperbolic group without non-trivial finite normal subgroups
(\cite[Ch. 5, Thm. 2.12]{Ghys}). Consequently, $E(\hat G)=\{1_{\hat G}\}$.

Now, consider the free product $F=\hat G * H$. $F$ is hyperbolic as a free product of hyperbolic groups
(\cite[Ch. 1, Exercise 4.34]{Ghys}) and non-elementary. Identify $\hat G$ and $H$ with their canonical copies inside
of $F$. Evidently, we have $E(\hat G)=E(F)=\{1_F\}$ in $F$, hence $\hat G$ is a $G$-subgroup of $F$.
By lemma \ref{non-periodic} one can find an element $g \in \hat G \le F$ of infinite order. Then
$$\langle g \rangle \cap H = \{1_F\}~\mbox{ in } F.$$

As it follows from the normal forms of elements of a free product, the subgroup $H$ is undistorted in $F$, hence, by
lemma \ref{quasi-undist}, $H$ is a quasiconvex subgroup of $F$. Define the quasiconvex subset
$Q \subset F$  by $Q=H \cup \langle g \rangle$. Obviously, no non-trivial element of $\hat G$
is conjugate to an element of $H$ in $F$, therefore, according to theorem \ref{qcex},
we can apply theorem \ref{mainthm} to obtain a non-elementary hyperbolic quotient $G_1$ of $F$ and an epimorphism
$\phi_0: F \to G_1$ that is surjective on $\hat G$, injective on $Q$, $\phi_0(H)$ is quasiconvex in $G_1$ and
\begin{equation} \label{E(G_1)} E(G_1)=\phi_0\bigl( E(F)\bigr)=\{1_{G_1}\}~,\end{equation}
\begin{equation} \label{<g>} \langle \phi_0(g) \rangle \cap \phi_0(H) = \{1_{G_1}\}~. \end{equation}
In particular, $\phi_0(H) \cong H~.$

Let the $\{\chi_j~|~j \in \N\}$ denote the set of all non-trivial conjugacy classes of elements in the group $G_1$.
Let $N_1$ be the normal subgroup of $G_1$ generated by $\chi_1$.
Observe that (\ref{E(G_1)}) implies that $N_1$ is infinite, consequently, it is non-elementary
(because $\Lambda(N_1)=\Lambda(G)=\partial G$ according to lemma \ref{limitofnormal} and this set is uncountable,
but the limit set of an infinite elementary subgroup consists of only two points).

By lemma \ref{elementalizereq} $E(N_1)=E(G_1)$ is trivial, hence, $N_1$ is a $G$-subgroup of the group $G_1$.
Denote $g_1 = \phi_0(g) \in G_1$, $H_1 =\phi_0(H) \le G_1$, $Q_1=\langle g_1 \rangle \cup H_1$.
The order of $g_1$ in the group $G_1$ is infinite,
hence (\ref{<g>}) implies that $|G_1:H_1|=\infty$. Therefore, $|N_1:(N_1 \cap h\langle g_1 \rangle h^{-1})|=\infty$
and $|N_1:(N_1 \cap hH_1h^{-1})|=\infty$ for any $h \in G_1$ (by lemma \ref{N-K}). Thus, by theorem \ref{qcex},
we can apply theorem \ref{mainthm} again and achieve a non-elementary hyperbolic quotient $G_2$ of $G_1$ together
with an epimorphism $\phi_1: G_1 \to G_2$ satisfying $\phi_1(N_1)=G_2$, $\phi_1$ is injective on $Q_1$,
$H_2=\phi_1(H_1)$ is a quasiconvex subgroup of $G_2$, $E(G_2)=\{1_{G_2}\}$ and
$\langle g_2 \rangle \cap H_2 = \{1_{G_2}\}$ where $g_2 = \phi_1(g_1)$.

Now, let $j_1=1$ and $j_2 > j_1$ be the smallest index such that $\phi_1(\chi_{j_2})$ is non-trivial in $G_2$.
Set $N_2 = \langle \phi_1(\chi_{j_2})\rangle \lhd G_2$. We can apply the same argument as before to get a
non-elementary hyperbolic quotient $G_3$ of $G_2$ with the natural epimorphism $\phi_2: G_2 \to G_3$
satisfying the properties we need (as above). And so on.

Thus, we obtain an infinite sequence of epimorphisms
$$G \stackrel{\psi}{\to} \hat G \stackrel{\phi_0}{\to} G_1 \stackrel{\phi_1}{\to} G_2 \stackrel{\phi_2}{\to} \dots~$$
where each epimorphism $\phi_{i}$ is injective on the image of $\phi_{i-1}(H)$, $i \in \N$.

Denote by $M$ the corresponding inductive limit of non-elementary hyperbolic groups. Then $M$ is a quotient of $G$.
As it is evident from the construction, $M$ is a simple group and the group $H$ is isomorphically embedded into $M$.
So, the corollary is proved. $\square$

\vspace{.15cm}
\underline{Proof of corollary \ref{monster}.} Let $A_1,A_2,A_3, \dots$ be an enumeration of all non-elemen\-tary
hyperbolic groups and $B_1,B_2,B_3,\dots$ -- an enumeration of all hyperbolic groups
(there are countably many of them since every hyperbolic group is finitely presented \cite{Mihalik}).
Denote $\A_i=A_i/E(A_i)$, $i=1,2,\dots$.

Set $F_1=\A_1 * B_1$.
Then, applying theorem \ref{mainthm}, we can obtain a non-elementary hyperbolic group $G_1$ and an epimorphism
$\phi_0: F \to G_1$ that is surjective on $\A_1$ and injective on $B_1$ (as before, we can demand that
$\phi_0(B_1)$ is quasiconvex in $G_1$, $|G_1:\phi_0(B_1)|=\infty$ and $E(G_1)=\{1_{G_1}\}$).

Again, let the $\{\chi_j~|~j \in \N\}$ be the set of all non-trivial conjugacy classes of elements in the group $G_1$,
$N_1 = \langle \chi_1 \rangle \lhd G_1$. By theorem \ref{mainthm} we obtain a (non-elementary hyperbolic) quotient
$\G_1$ with the natural epimorphism $\psi_1: G_1\to \G_1$ that is surjective on $N_1$ and injective on the image
$\phi_0(B_1)$ of $B_1$ in $G_1$.

Next, define $F_2=\G_1*\A_2*B_2$. Let $G_2$ be a non-elementary hyperbolic quotient of $F_2$ such that the natural
epimorphism $\phi_1:F_2 \to G_2$ is surjective on the subgroups $\G_1,\A_2 \le F_2$ and injective on $B_2$ and
the image of $B_1\le G_1\le F_2$.

Now, let $j_1=1$ and $j_2 > j_1$ be the smallest index such that the image of $\chi_{j_2}$ (under the composition
$\phi_1\circ \psi_1$) is non-trivial in $G_2$. Let $N_2 = \langle (\phi_1\circ \psi_1)(\chi_{j_2})\rangle \lhd G_2$.
Then we can find an epimorphism $\psi_2: G_2 \to \G_2$ onto a non-elementary hyperbolic group $\G_2$
that is surjective on $N_2$ and injective on the images of $B_1,B_2$.

And so on. Thus we achieve a sequence of epimorphisms
$$\hat G_1 \to \hat G_2 \to \G_3 \to \dots~.$$
Let $M$ the corresponding inductive limit of these groups. As it follows from the construction, $M$ satisfies
all the properties required. $\square$

\section{Thrifty Embeddings}
Observe that since any elementary torsion-free group is cyclic, maximal elementary subgroups are malnormal in
torsion-free non-elementary hyperbolic groups.

\vspace{.15cm}
\remark \label{malnrem} Let $H$ be a malnormal subgroup of a group $G$, $g \in G$. Then

$(a)$ The conjugate subgroup $gHg^{-1} \le G$ is also malnormal;

$(b)$ If $K \le G$ is an infinite subgroup and $|K:(K\cap gHg^{-1})|<\infty$ then $g^{-1}Kg \le H$;

$(c)$ For any $h \in H\backslash\{1_G\}$, $C_G(h) \le H$;

$(d)$ If $f \in G$ and $fHf^{-1} \cap gHg^{-1}\neq \{1_G\}$ then $fHf^{-1}=gHg^{-1}$.

\vspace{.15cm}

For the proof of corollary \ref{thrifty} we will need the following auxiliary lemma:

\begin{lemma} \label{HNN-maln} Suppose $G$ is a group and $H,A,B$ are its subgroups.
Assume that $H$ and $B$ are malnormal in $G$, $H \cap gBg^{-1}=\{1_G\}$ for any $g \in G$ and
there is an isomorphism $\tau: A \to B$. Then the natural image of $H$ in the HNN-extension
$$G_1=\langle G,t~|~tAt^{-1}=B\rangle \stackrel{def}{=}\langle G,t~|~tat^{-1}=\tau(a),~a \in A\rangle$$
is malnormal.
\end{lemma}

\underline{Proof.}
%
Identify $G$ and $H$ with their canonical images in $G_1$. Assume that there exists
$w \in G_1\backslash H$ and non-trivial elements $x,y \in H$ such that $wxw^{-1}=y$.
Then we can write \begin{equation} \label{w}
w=u_0t^{\epsilon_1}u_1t^{\epsilon_2}\cdot \dots \cdot t^{\epsilon_{n-1}}u_{n-1}t^{\epsilon_n}u_n~~\mbox{ in $G_1$ }~,
\end{equation}
where $u_0,u_n \in G$, $u_1,\dots,u_{n-1} \in G\backslash\{1_G\}$, $\epsilon_1,\dots,\epsilon_n \in \{1,-1\}$,
and this representation is reduced (i.e. it contains no occurences of the form
$tut^{-1}$ or $t^{-1}vt$ where $u \in A$, $v \in B$).


Observe that $n \ge 1$ since $w \notin G$ (by malnormality of $H$ in $G$) and
\begin{equation} \label{u_nx}   u_0t^{\epsilon_1}\cdot \dots \cdot t^{\epsilon_n}u_n x
u_n^{-1}t^{-\epsilon_n}\cdot \dots\cdot t^{-\epsilon_2}u_0^{-1} y^{-1}=1_{G_1}~.\end{equation}


By Britton's lemma (\cite{L-S}) the left-hand side in (\ref{u_nx}) is not reduced, hence
$u_n x u_n^{-1}$ belongs to $A$ or $B$. But this element is a conjugate of $x\in H$ therefore, according to the
assumptions of the lemma, it has to be in $A$ and $\epsilon_n=1$. Consequently,
$t^{\epsilon_n} u_n x u_n^{-1}t^{-\epsilon_n}=v\in B \backslash\{1_G\}$. 
Since no element of $B$ is conjugate to the element $y \in H$
in the group $G$, the number $n$ from the representation (\ref{w}) must be at least $2$ and
$$w  x w^{-1} y^{-1} \stackrel{G_1}{=} u_0t^{\epsilon_1}\cdot \dots \cdot t^{\epsilon_{n-1}}u_{n-1} v
u_{n-1}^{-1}t^{-\epsilon_{n-1}}\cdot \dots\cdot t^{-\epsilon_2}u_0^{-1} y^{-1}=1_{G_1}~.$$

Applying Britton's lemma again, we get that the element $u_{n-1} v u_{n-1}^{-1}$ either belongs to
$A$ (and $\epsilon_{n-1}=1$) or to $B$ (and $\epsilon_{n-1}=-1$).
So, if it is in $A$, then $t^{\epsilon_{n-1}}u_{n-1} v u_{n-1}^{-1}t^{-\epsilon_{n-1}} \in B$ and
$n$ has to be at least $3$; thus we can proceed as before. This process will end after finitely many steps
because each time we eliminate a $t^{\pm 1}$-element from the representation (\ref{w}) of $w$.
Therefore, we can assume that
$u_{n-1} v u_{n-1}^{-1}\in B$ and $\epsilon_{n-1}=-1$. But the subgroup $B$ was malnormal in
$G$ and $v \in B\backslash\{1_G\}$, hence $u_{n-1}\in B$. Hence
$t^{\epsilon_{n-1}}u_{n-1}t^{\epsilon_n}\equiv t^{-1}u_{n-1}t \in A$ which contradicts to our assumption
that the right-hand side of (\ref{w}) is reduced.
%
%

The lemma is proved. $\square$

\begin{lemma} {\normalfont{(\cite[Thm. 3]{Mihaj-Olsh},\cite[Cor. 1]{Khar-Myas})}} \label{HNN-hyp}
Let $G$ be a hyperbolic group with isomorphic infinite elementary
subgroups $A$ and $B$, and let $\tau$ be an isomorphism from $A$ to $B$. The HNN-extension
$G_1=\langle G,t~|~tat^{-1}=\tau(a),a \in A\rangle$ of $G$ with associated subgroups $A$ and $B$ is hyperbolic if
and only if the following two conditions hold:\\
$1)$ either $A$ or $B$ is a maximal elementary subgroup of $G$;\\
$2)$ for all $g\in G$ the subgroup $gAg^{-1} \cap B$ is finite.
\end{lemma}

\begin{lemma} {\normalfont{(\cite[Thm. 4]{Khar-Myas})}} \label{HNN-qc} Let the HNN-extension
$G_1=\langle G,t~|~tAt^{-1}=B\rangle$ be hyperbolic with $A$ quasiconvex in $G_1$. Then $G$ is quasiconvex in $G_1$.
\end{lemma}

\remark \label{qc_in_qc} Suppose $G_1$ is a hyperbolic group and $H \le G \le G_1$. If $H$ is quasiconvex in $G$ and
$G$ is quasiconvex in $G_1$ then $H$ is quasiconvex in $G_1$.

\vspace{.15cm}
This follows from lemma \ref{quasi-undist} and the observation that an undistorted subgroup of an undistorted subgroup
is undistorted in the entire group.

We are now ready to give the

\vspace{.15cm}
\underline{Proof of corollary \ref{thrifty}.} Consider the free product $F=G*H$. Then $F$ is a non-elementary torsion-free
hyperbolic group, $G$ is a $G$-subgroup of $F$ and $H$ is quasiconvex in $F$ (because it is undistorted).
H is non-trivial by the assumptions of the corollary, hence there is an element $y \in H$ of infinite order.
Pick any $f \in G\backslash \{1_F\}$ and set $x=fyf^{-1}\in F$. From normal
forms of elements of the free product $F$ it follows that $H$ is malnormal in $F$,
$gHg^{-1}\cap G=\{1_F\}$ for any $g \in F$ and the infinite cyclic subgroup of $F$ generated by $x$
has trivial intersection with $H$. Denote $Q=\langle x \rangle \cup H$ -- a quasiconvex subset of $F$.

By theorems \ref{qcex} and \ref{mainthm} there exists a non-elementary hyperbolic quotient $G_0$
of $F$ and an epimorphism $\psi_0:F \to G_0$ with the properties $1)-9)$ from the claim of theorem \ref{mainthm}.
Thus $\psi_0(G)=G_0$, $\psi_0$ is injective on $Q$, $G_0$ is torsion-free (by the property $7)$),
$\psi_0(H)$ is quasiconvex in $G_0$, $\psi_0(x) \in (G_0)^0$ and $\psi_0(H) \cap \langle \psi_0(x) \rangle =\{1_{G_0}\}$.

Suppose for some non-trivial $z \in G_0$ there are non-trivial $a,b \in H$
such that $z\psi_0(a)z^{-1}=\psi_0(b)$. By property $4)$ from the claim of theorem \ref{mainthm}, there exists
an element $u \in F$ such that $uau^{-1}=b$. $H$ was malnormal in $F$, therefore $u \in H$ and
$z^{-1}\psi_0(u)\psi_0(a)\psi_0(u)^{-1}z=\psi_0(a)$, i.e. $z^{-1}\psi_0(u) \in C_{G_0}\bigl(\psi_0(a)\bigr)$.
Then, according to property $5)$, there is $v \in C_G(a)$ satisfying $\psi_0(v)=z^{-1}\psi_0(u)$. Also, by
remark \ref{malnrem}, $v \in H$. Thus, $z = \psi_0(u)\psi_0(v)^{-1} \in \psi_0(H)$, i.e. $\psi_0(H)$ is malnormal
in $G_0$.

Enumerate all non-trivial elements of the group $G_0$: $g_1,g_2,\dots$,
and all its two-generated non-elementary subgroups: $K_1,K_2,\dots$.

The group $M$ will be constructed as an inductive limit of groups $G_i$, $i=0,1,\dots$.
Assume, the non-elementary hyperbolic torsion-free  quotient  $G_{i-1}$ of $G_0$ has already been constructed,
$i \ge 1$, and it satisfies the following properties: the natural epimorphism
$\pi_{i-1}:G_0 \to G_{i-1}$ ($\pi_0=id_{G_0}:G_0 \to G_0$) is injective on $\psi_0(H) \cup \langle \psi_0(x) \rangle$,
the image of $\psi_0(H)$ is quasiconvex and malnormal in $G_{i-1}$; images of the elements $g_1,\dots,g_{i-1}$
are conjugate in $G_{i-1}$ to some elements from $\pi_{i-1}\bigl(\psi_0(H)\bigr)$, and images of the subgroups
$K_1,\dots,K_{i-1}$ either coincide with $G_{i-1}$ or are conjugate in $G_{i-1}$ to a subgroup of
$\pi_{i-1}\bigl(\psi_0(H)\bigr)$, or are elementary.

Let us now construct the group $G_i$. Consider the element $\pi_{i-1}(g_i) \in G_{i-1}$.
To simplify the notation, identify $H$ and $\pi_{i-1}\bigl(\psi_0(H)\bigr)$.
If $\pi_{i-1}(g_i)$ is conjugate in $G_{i-1}$ to an element from $H$ then set $F_i=G_{i-1}$.

If not, then the element $\pi_{i-1}(g_i)$ has infinite order in $G_{i-1}$ and the maximal elementary subgroup
$B=E\bigl(\pi_{i-1}(g_i)\bigr)$ is infinite cyclic (because $G_{i-1}$ is torsion-free)
and malnormal in $G_{i-1}$. Part $(b)$ of remark
\ref{malnrem} implies that $H\cap gBg^{-1}=\{1_{G_{i-1}}\}$ for any $g \in G_{i-1}$.
Denote by $A \le G_{i-1}$ the infinite cyclic subgroup
of $H$ generated by the element $y$ chosen in the beginning of the proof.
Then we can construct an HNN-extension $$F_i=\langle G_{i-1},t~|~tAt^{-1}=B\rangle~. $$
According to lemmas \ref{HNN-hyp} and \ref{HNN-qc} and basic properties of HNN-extensions (see \cite[Ch. IV]{L-S}),
$F_i$ is a torsion-free non-elementary hyperbolic group and the natural image of $G_{i-1}$ is quasiconvex in it.
By lemma \ref{HNN-maln} and remark \ref{qc_in_qc}, $H$ is malnormal and quasiconvex in $F_i$.
Note that the latter implies that $|G_{i-1}:(G_{i-1} \cap gHg^{-1})|=\infty$ for any $g \in F_i$ because, otherwise,
by part $(b)$ of remark \ref{malnrem}, $G_{i-1} \le gHg^{-1}$ and since $H\le G_{i-1}$ is non-trivial, part
$(d)$ of the same remark would claim that $G_{i-1}\le gHg^{-1}=H$. This leads to a contradiction with the fact that
$x \in G_{i-1} \backslash H$.

By construction, $\pi_{i-1}(g_i)$ is conjugate to some element of $H$ in $F_i$.

Now consider the subgroup $\pi_{i-1}(K_i)\le G_{i-1} \le F_i$. If this subgroup is elementary or conjugate to a subgroup
of $H$ in $F_i$, then we apply theorems \ref{qcex} and \ref{mainthm} to obtain a torsion-free non-elementary hyperbolic
group $G_i$ and an epimorphism $\psi_i:F_i \to G_i$ such that $\psi_i(G_{i-1})=G_i$, $\psi_i$ is injective on
$H \cup \langle x \rangle$, $\psi_i(H)$ is quasiconvex and malnormal (as before) in $G_i$. Then
$\psi_i\bigl(\pi_{i-1}(K_i)\bigr) \le G_i$ is either elementary or conjugate to a subgroup of $\psi_i(H)$ in $G_i$.

Thus, we can assume that $\pi_{i-1}(K_i)$ is non-elementary and not conjugate with a subgroup of $H$ in $F_i$.
Then, by remark \ref{malnrem}, $$|\pi_{i-1}(K_i):(\pi_{i-1}(K_i) \cap gHg^{-1})|=\infty \quad
\mbox{for any $g \in F_i$},$$  hence
we can use theorems \ref{qcex} and \ref{mainthm} to get an epimorphism $\psi_i$ of $F_i$ onto a
non-elementary torsion-free hyperbolic group $G_i$ satisfying the following conditions:
$\psi_i\bigl(\pi_{i-1}(K_i)\bigr)=G_i$ (consequently, $\psi_i(G_{i-1})=G_i$), $\psi_i$ is injective on
$H \cup \langle x \rangle$ and $\psi_i(H)$ is quasiconvex and malnormal in $G_i$.

Thus, we have constructed the group $G_i$ for every $i=0,1,2,\dots$.

Set $M=\lim_{\to}(G_i,\psi_{i+1})$. It remains to prove that $M$ satisfies
the properties required. There is a natural epimorphism $\pi:G_0\to M$.
Note that if a word $w$ is trivial in $M$, then (by the definition of an inductive limit)
$w$ is trivial in $G_i$ for some $i$, hence $M$ is torsion-free, $\pi$ is injective on $H$, $\pi(x)\neq 1_M$ and
$\pi(H) \cap \langle \pi(x) \rangle=\{1_M\}$ (we identify $H$ and $x$ with their images in $G_0$). Therefore
$\pi(H)$ is a proper subgroup of $M$ and, since the image of $H$ was malnormal in each $G_i$,
$\pi(H)$ will be malnormal in $M$.

Denote $P=\pi(H) \le M$ and assume that $L$ is a proper non-trivial subgroup of $M$.
Then there exists $a \in L\backslash \{1_M\}$. Suppose that for every $b \in L$ there exists
$g_b \in M$ such that $g_b\langle a,b \rangle g_b^{-1} \le P$. Set $g=g_{1_M}$ and pick an arbitrary $b \in L$.
Then $$a \in g_b^{-1}Pg_b \cap g^{-1}Pg \neq \{1_M\}~.$$ Therefore, applying remark \ref{malnrem}, we get
$g_b^{-1}Pg_b = g^{-1}Pg$, thus $b \in g^{-1}Pg$ for any $b \in L$, hence $gLg^{-1} \le P$.
So, if $L$ is not conjugate to a subgroup of $P$ then there should exist $b \in L$ such that the subgroup
$\langle a,b \rangle \le L$ is not conjugate to any subgroup from $P$.
Choose arbitrary elements $c,d \in G_0$ with $\pi(c)=a$, $\pi(d)=b$. Then
$\pi(\langle c,d \rangle)=\langle a,b \rangle$ and the image of $\langle c,d \rangle$ in $G_i$ is not conjugate to
a subgroup of the (corresponding) image of $H$ for all $i$. Thus, this image is non-elementary (i.e. non-cyclic) in $G_i$
for all $i$ (since every cyclic subgroup will eventually be conjugate to some cyclic subgroup from an image of $H$).
Consequently, $\langle c,d \rangle=K_j$ for some $j \in \N$ and the homomorphism $\pi_j:G_0\to G_j$ will be
surjective on $K_j$. It follows that $\langle a,b \rangle=\pi(\langle c,d \rangle)=M$ -- a contradiction
with the condition $L\neq M$. So, we showed that any proper subgroup $L$ of $M$ is conjugate to some subgroup
of $P$.

Finally, if $N \lhd M$ and $N\neq M$ then, applying the above, we obtain an element $g \in M$ such that
$N=gNg^{-1}\le P$. But this implies that $N=\{1_M\}$ because $P$ is malnormal. Thus, $M$ is simple. $\square$

\end{document}